%% file: main.tex
\documentclass[12pt]{article}
\usepackage[a4paper, left=3cm, right=3cm, top=3cm, bottom=3cm]{geometry}
\usepackage{amssymb, amsmath, amsthm, stmaryrd, enumerate}
\usepackage[ps,dvips,all]{xy}
\SelectTips{cm}{12}

\DeclareMathOperator{\ob}{ob}
\newcommand{\db}[1]{\mathbb{#1}}
\newcommand{\cat}[1]{\mathbf{#1}}
\newcommand{\op}{\mathrm{op}}

\newcommand{\id}{\mathrm{id}}
\newcommand{\thg}{{\mathord{\text{--}}}}

\newcommand{\cotimes}{\mathbin{\vphantom{\times}_s\mathord{\times}_{t}}}

\newcommand{\spn}[1]{{\left<{#1}\right>}}
\newcommand{\elt}[1]{\left\ulcorner{#1}\right\urcorner}

\newcommand{\defn}[1]{\textbf{#1}}
\newcommand{\cd}[2][]{\vcenter{\hbox{\xymatrix#1{#2}}}}
\newcommand{\cdl}[2][]{\xymatrix@1#1{#2}}

\renewcommand{\b}[1]{\boldsymbol{\mathbf{#1}}}
\newcommand{\I}{{\b I}}

\newcommand{\A}{{\b A}}
\newcommand{\B}{{\b B}}
\newcommand{\C}{{\b C}}
\newcommand{\D}{{\b D}}

\newcommand{\W}{{\b W}}
\newcommand{\X}{{\b X}}
\newcommand{\Y}{{\b Y}}
\newcommand{\Z}{{\b Z}}
\renewcommand{\a}{{\mathfrak{a}}}
\newcommand{\e}{{\mathfrak{e}}}
\renewcommand{\l}{{\mathfrak{l}}}
\newcommand{\m}{{\mathfrak{m}}}
\renewcommand{\r}{{\mathfrak{r}}}

\newcommand{\xtor}[1]{\cdl[@1]{{} \ar[r]|-{\object@{|}}^{#1} & {}}}
\newcommand{\tor}{\ensuremath{\relbar\joinrel\mapstochar\joinrel\rightarrow}}
\newcommand{\To}{\ensuremath{\Rightarrow}}
\newcommand{\Tor}{\ensuremath{\Relbar\joinrel\Mapstochar\joinrel\Rightarrow}}

\newcommand{\twocong}[2][0.5]{\ar@{}[#2] \save ?(#1)*{\cong}\restore}
\newcommand{\rtwocell}[3][0.5]{\ar@{}[#2] \ar@{=>}?(#1)+/l 0.2cm/;?(#1)+/r 0.2cm/^{#3}}
\newcommand{\ltwocell}[3][0.5]{\ar@{}[#2] \ar@{=>}?(#1)+/r 0.2cm/;?(#1)+/l 0.2cm/_{#3}}
\newcommand{\dtwocell}[3][0.5]{\ar@{}[#2] \ar@{=>}?(#1)+/u  0.2cm/;?(#1)+/d 0.2cm/^{#3}}
\newcommand{\dthreecell}[3][0.5]{\ar@{}[#2] \ar@3{->}?(#1)+/u  0.2cm/;?(#1)+/d 0.2cm/^{#3}}
\newcommand{\utwocell}[3][0.5]{\ar@{}[#2] \ar@{=>}?(#1)+/d 0.2cm/;?(#1)+/u 0.2cm/_{#3}}
\newcommand{\dtwocelltarg}[3][0.5]{\ar@{}#2 \ar@{=>}?(#1)+/u  0.2cm/;?(#1)+/d 0.2cm/^{#3}}
\newcommand{\utwocelltarg}[3][0.5]{\ar@{}#2 \ar@{=>}?(#1)+/d  0.2cm/;?(#1)+/u 0.2cm/_{#3}}

\newtheorem{Prop}{Proposition}
\newtheorem{Cor}[Prop]{Corollary}

\theoremstyle{definition}
\newtheorem{Defn}[Prop]{Definition}

\newtheorem{Ex}[Prop]{Example}

\newcommand{\coll}{{\db Coll}}
\newcommand{\ten}{\mathbin{\bullet}}
\newcommand{\bten}{\mathbin{\bullet}}

\newcommand{\unit}{e}
\newcommand{\dCat}{{\db Cat}}
\newcommand{\Cat}{{\cat{Cat}}}
\newcommand{\ks}{\db L / S\I_1}

\makeatletter \@namedef{itemize*}{\itemize\parsep\z@ \parskip\z@}
\@namedef{enditemize*}{\enditemize} \@namedef{enumerate*}{\enumerate\parsep\z@
\parskip\z@} \@namedef{endenumerate*}{\endenumerate} \makeatother

\begin{document}

\title{Double clubs}
\author{Richard Garner\\
\fontsize{12}{14}\selectfont\emph{University of Cambridge,}\\
\fontsize{12}{14}\selectfont\emph{Department of Pure Mathematics and Mathematical Statistics,}\\
\fontsize{12}{14}\selectfont\emph{Cambridge CB3 0WB, UK}\\
\fontsize{12}{14}\selectfont\texttt{rhgg2@cam.ac.uk}}
\maketitle

\input abstract.tex

\input intro.tex
\input chapter1.tex

\input chapter2.tex

 \input chapter3.tex
\input chapter4.tex

\input chapter5.tex

\input chapter6.tex
\input appendix.tex
\bibliography{biblio}
\end{document}

%% file: abstract.tex
\noindent \textbf{Abstract}. We develop a theory of \emph{double clubs} which extends Kelly's
theory of clubs to the pseudo double categories of Par\'e and
Grandis. We then show that the club for symmetric strict monoidal
categories on $\cat{Cat}$ extends to a `double club' on the pseudo
double category $\dCat$ of `categories, functors, profunctors and
transformations'.

%% file: intro.tex
\section{Introduction}

Kelly's theory of \emph{clubs} \cite{kelly:clubsold2, kelly:clubsold1,
  kelly:clubsold3, kelly:clubs} captures an important intuition, that
of adding structure in a `generic way'. In the case of $\cat{Cat}$, it
tells us that, given a description of this added structure at the
terminal category $1$, we should be able to derive it at an arbitrary
category $\cat C$ by `labelling with objects and maps of $\cat C$'.

The genesis of this paper was an attempt to do something similar for
$\cat{Mod}$, the bicategory of categories and profunctors. As it
stands, the theory of clubs is inadequate: it deals with categories
with pullbacks, whilst $\cat{Mod}$ is neither a category nor has pullbacks.
Therefore, we must look for a suitable generalisation of the
theory of clubs which is amenable to application in $\cat{Mod}$.

Now, taking pullbacks is fundamental to the theory of
clubs, so we are led to question whether or not $\cat{Mod}$ is the
correct place to work; ideally, we should like to replace it with
something where we \emph{can} take lots of pullbacks.
Now, observe that $\cat{Mod}$ has certain peculiar properties: it
has all lax colimits, but these lax colimits have a universal
property up to \emph{isomorphism} rather than up to
\emph{equivalence}; unfortunately, the language of bicategories
cannot express what this universal property is. Similarly, the
operation given on objects by cartesian product of categories
induces a structure of monoidal bicategory on $\cat{Mod}$; again,
this structure ought to be associative up to \emph{isomorphism}
rather than \emph{equivalence}, and again, the language of
bicategories is simply unable to express this.

Inspired by this, we are led to consider the \emph{pseudo double
categories} of \cite{grandispare:limits} and
\cite{grandispare:adjoints} (and also considered briefly by
\cite{leinster:higheroperads}). These are a weakening of
Ehresmann's notion of \emph{double category}
\cite{ehresmann:catstruct, ehresmann:catbook}, and have two
directions, one `category-like' and the other `bicategory-like'.
The presence of a `category-like' direction allows us to express
`up-to-isomorphism' as well as `up-to-equivalence' notions, and
more saliently, to take lots of pullbacks. Indeed, in our case, we
can generalise $\cat{Mod}$ to the pseudo double category $\dCat$
of `categories, functors, profunctors and transformations' which
in an appropriate sense, has \emph{all} pullbacks.

The main thrust of this paper, then, is to develop a suitable
generalisation of the theory of clubs from plain categories to pseudo
double categories. Concurrently, we generalise the leading example of
a club on $\cat{Cat}$, the club for symmetric strict monoidal
categories, to such a `double club' on the pseudo double category
$\dCat$.

This paper is not mere theory for theory's sake:
it has been developed very much with an application in mind.  In
\cite{polycats}, we make extensive use of these results to get a
handle on the ``higher-dimensional bookwork'' involved in the
construction of a \emph{pseudo-distributive law} \cite{tanaka:thesis}
on $\cat{Mod}$. An examination of \cite{polycats}, therefore, may
give the reader a better feel for the motivation behind the present work.

\textbf{Structure of the paper}.  In Section 2, we summarise the basic
concepts and definitions of pseudo double categories, and prove some
new results about double functor categories. In Section 3, we recap
the theory of plain clubs, before, in Section 4, starting our
generalisation of this theory to the setting of double categories.
First, we explore some necessary further aspects of the theory of
pseudo double categories, considering slice double categories,
equivalences of double categories and cartesian maps in double
categories, and then prove a key equivalence of double categories.

In Section 5, we develop the theory of `monoidal
double categories': with this in place, we are ready, in Section 6, to
give two definitions of `double club', one more abstract, the
other more tractable. Finally, in Section 7, we show that we can
extend the club $S$ for symmetric strict monoidal categories on
$\cat{Cat}$ to a double club on $\dCat$.

Two Appendices gives a result on equivalences in double categories
(Appendix A), and a technical result on `whiskering' which is of some
use in applying the theory of double clubs (Appendix B).

\textbf{Acknowledgements}. I should like to thank Martin Hyland for illuminating discussions, and
an anonymous referee for perceptive and helpful remarks on a earlier
draft of this paper. This work was supported by a PhD grant from the
EPSRC.


%% file: chapter1.tex
\section{Pseudo double categories}
We begin by recapping some of the theory of \emph{pseudo double
  categories}. Since the full details of this can be found in
\cite{grandispare:limits, grandispare:adjoints}, we shall
merely set out our notation and give a few examples.

\subsection{Basic theory}
\begin{Defn}
A \defn{pseudo double category} $\db K$ consists of:
\begin{itemize*}
\item A diagram of categories \smash{$\xymatrix@1{K_1 \ar@<1.5pt>[r]^s \ar@<-1.5pt>[r]_t & K_0}$}. We write $X$ for a typical object and $f$ for a typical arrow of
$K_0$, and call them \defn{objects} and
\defn{vertical maps} of $\db K$; similarly, we write $\X$ for a typical object
and $\b f$ for a typical arrow of $K_1$, and call them \defn{horizontal maps}
and \defn{cells} of $\db K$. We call $s$ and $t$ the \defn{source} and
\defn{target} functors of $\db K$, and write $X_s$ and $X_t$ for $s(\X)$ and
$t(\X)$; similarly, we write $f_s$ and $f_t$ for $s(\b f)$ and $t(\b f)$.
\item A \defn{horizontal units} functor
$\I \colon K_0 \to K_1$.
\item A \defn{horizontal composition} functor
$\otimes \colon K_1
\mathbin{\vphantom{\times}_s\mathord{\times}_{t}} K_1 \to K_1$,
where $K_1
\mathbin{\vphantom{\times}_s\mathord{\times}_{t}} K_1$ is the evident
pullback.
\item Special isomorphisms
\begin{gather*}\l_{\X} \colon \X \to \I_{X_t} \otimes \X\text,\quad  \r_{\X} \colon \X \to \X
\otimes \I_{X_s}\text,\\
 \text{and} \quad
\a_{\X\Y\Z} \colon \X \otimes (\Y \otimes \Z) \to (\X \otimes \Y) \otimes \Z\text,\end{gather*}
in $K_1$, natural in all variables.
Note that we say a map $\b f \colon \X \to \Y$ in $K_1$ is
\defn{special} if $X_s = Y_s$, $X_t = Y_t$, $f_s = \id_{X_s}$ and $f_t =
\id_{X_t}$.
\end{itemize*}
These data are required to satisfy five straightforward axioms.
\end{Defn}
\noindent We write a typical object of $\db K$ as $X$, a typical vertical arrow
as $f \colon X \to Y$, a typical horizontal arrow as $\X \colon X_s \tor X_t$
and a typical cell as
\[\cd{
X_s \ar[d]_{f_s} \ar[r]|-{\object@{|}}^{\X} \dtwocell{dr}{\b f}& X_t \ar[d]^{f_t}  \\
Y_s \ar[r]|-{\object@{|}}_{\Y} & Y_t\text, }\]
which we may abbreviate to $\b f \colon \b X \Rightarrow \b Y$.
We observe that any pseudo double category $\db K$ contains a
bicategory $\mathcal B \db K$, with objects the objects of $K_0$, 1-cells the
objects of $K_1$ and 2-cells the special maps in $K_1$. Therefore, given $\X \colon
A \tor B$ and $\Y \colon B \tor C$ in $K_1$, we may notate $\Y \otimes \X$ as
\[\cdl{{\Y \otimes \X} \colon A \ar[r]|-{\object@{|}}^-{\X} & B \ar[r]|-{\object@{|}}^-{\Y} & C\text,}\]
and can extend this notation to horizontal composition of cells.  As
for bicategories, this notation is ambiguous for chains of three or
more such composites: any such will need a choice of bracketing in
order to specify a composite horizontal arrow of $\db K$.  However, as
for bicategories, we may use pasting diagrams to specify composites of
special maps in $K_1$: it follows from the bicategorical pasting
theorem \cite{power:pasting, verity:thesis} that such diagrams
uniquely specify a special map in $K_1$ once a bracketing for the
start and end edge has been chosen.

\begin{Ex}\label{dcat}
The pseudo double category $\dCat$ is given as follows:
\begin{itemize*}
\item \textbf{Objects} are small categories $X$, $Y$, \dots;
\item \textbf{Vertical maps} are functors $F \colon X \to Y$;
\item \textbf{Horizontal maps} $\X \colon X_s \tor X_t$ are
profunctors from $X_s$ to $X_t$; i.e., functors $\X \colon X_t^\op \times X_s \to \cat{Set}$. We shall specify such by giving:
\begin{itemize*}
\item The \textbf{proarrows} $g \colon x_t \tor x_s$, for $x_t \in X_t$ and $x_s \in X_s$: in other words, the elements of $\X(x_t; x_s)$;
\item The \textbf{actions} by maps of $X_t$ and $X_s$; so for $h \colon
  x_s \to x_s'$ in $X_s$ and $f \colon x_t' \to x_t$ in $X_t$,
  we give the functors
\[
h \bullet (\thg) = \X(\id_{x_t}; h) \colon \X(x_t; x_s) \to \X(x_t; x_s')\] and \[(\thg)
\bullet f = \X(f; \id_{x_s}) \colon \X(x_t; x_s) \to \X(x_t'; x_s)\text.
\]
Given a proarrow $g \colon x_t \tor x_s$, we write the
elements $h \bullet g$ and $g \bullet f$ as
\[
\cdl{x_t \ar[r]|-{\object@{|}}^{g} & x_s \ar[r]^h & x_s'} \quad \text{ and } \quad
\cdl{x_t' \ar[r]^{f} & x_t \ar[r]|-{\object@{|}}^g & x_s}
\]
respectively. By analogy with categorical composition, we'll tend to
drop the `$\bullet$' symbol where convenient, and denote these
actions simply by juxtaposition;
\end{itemize*}
\item \textbf{Cells} $\b F \colon \b X \Rightarrow \b Y$
are natural transformations
\[\cd{
 X_t^\op \times X_s
  \ar[rr]^{F_t^\op \times F_s}
  \ar[dr]_{\X} &
 {}
  \rtwocell{d}{\b F} &
 Y_t^\op \times Y_s
  \ar[dl]^{\Y} \\ &
 \cat{Set}\text.
}\]
We shall specify a cell by giving its action on proarrows of $\b X$; in other words, by giving the components
\[{\b F}_{x_t, x_s} \colon \X(x_t; x_s) \rightarrow \Y\big(F_t(x_t); F_s(x_s)\big)\text.\]
In practice, we drop the suffices and refer to all of
these maps simply as `$\b F$'. Note that naturality of $\b F$ amounts to
verifying the equivariance formulae
\[\b F(h \bullet g) = F_s(h) \bullet \b F(g) \quad \text{and} \quad
\b F(g \bullet f) = \b F(g) \bullet F_t(f)\text.\]
\end{itemize*}
Vertical composition is given as in $\cat{Cat}$, whilst horizontal
composition $\otimes$, horizontal units $\I$, associativity $\a$ and
unitality $\l$, $\r$ are given as in $\cat{Mod}$, the bicategory of
categories and profunctors. In particular, we notate the proarrows
of $\I_X$ (the identity at $X$) by
\[\I_f \colon x \tor y \quad \text{where} \quad f \in \X(x, y)\]
and the proarrows of $\cdl{{\Y \otimes \X} \colon A
  \ar[r]|-{\object@{|}}^-{\X} & B \ar[r]|-{\object@{|}}^-{\Y} & C}$ by
\[k \otimes g \colon c \tor a \quad \text{where} \quad k \in \Y(c; b), \quad
g \in \X(b; a)\text.\] Note that in the latter case, the `proarrows'
are subject to the equivalence relations $gf \otimes k \simeq g
\otimes fk$ for suitable $f \in B(b, b')$; as usual we shall conflate
$k \otimes g$ with its image under this equivalence relation.

From this example, we can derive several more useful examples: we can
restrict our attention to the \emph{discrete} categories, to get the
pseudo double category of sets, maps and spans; we can replace
categories with $\mathcal V$-categories (for some suitable base for enrichment
$\mathcal V$) to produce the pseudo double category $\mathcal V$-$\dCat$; and we can
restrict this last to \emph{one-object} $\mathcal V$-categories, thereby
producing the pseudo double category of monoids, monoid maps and
modules in $\mathcal V$. In particular, setting $\mathcal V = \cat{Ab}$, the category
of abelian groups, we get the pseudo double category of rings, ring
homomorphisms and bimodules.
\end{Ex}

\begin{Defn}
A \defn{morphism of pseudo double categories} (or \defn{double
morphism} for short) $F \colon \db K \to \db L$ consists of
functors $F_0 \colon K_0
\to L_0$ and $F_1 \colon K_1 \to L_1$~--~and to ease notation we write
`$F$' interchangeably for both~--~together with
special maps
$\e_X \colon \I_{FX} \to F\I_X$ and  $\m_{\X, \Y} \colon F\X \otimes F\Y \to F(\X \otimes
\Y)$, natural in all variables, all satisfying five evident axioms.
\end{Defn}

Pseudo double categories and the morphisms between them form
themselves into a category $\cat{DblCat}$. Similarly, we may define
the category $\cat{DblCat}_o$ of `pseudo double categories and double
opmorphisms' and $\cat{DblCat}_\psi$ of `pseudo double categories and
homomorphisms': for an opmorphism, $\e_X$ and $\m_{\X, \Y}$ point in
the opposite direction, whilst for a homomorphism, $\e_X$ and $\m_{\X,
  \Y}$ are invertible.

\begin{Ex}\label{S}
  We give an example of a homomorphism on the pseudo double category
  $\dCat$ of Example \ref{dcat}.
  This homomorphism $S \colon \dCat \to \dCat$ extends the
  `free symmetric strict monoidal category 2-functor' on $\cat{Cat}$,
  and given as follows:
\begin{itemize}
\item \textbf{On objects}: Given a small category $X$, the category $SX$ has:
\begin{itemize}
\item \textbf{Objects} being pairs $(n, \spn{x_i})$, where $n \in \mathbb N$ and $x_1, \dots,
x_n \in \ob X$;
\item \textbf{Arrows} being
\[(\sigma, \spn{g_i}) \colon (n, \spn{x_i}) \to (n, \spn{y_i})\text,\] where
$\sigma \in S_n$ and $g_i \colon x_i \to y_{\sigma(i)}$ in $X$ (note
that there are no maps from $(n, \spn{x_i})$ to $(m, \langle{y_j}\rangle)$ for $n \neq m$).
\end{itemize}
Composition and identities in $SX$ are given in the evident
way; namely,
\begin{align*}
\id_{(n, \spn{x_i})} &= (\id_n, \spn{\id_{x_i}})\\
\text{and } (\tau, \spn{g_i}) \circ (\sigma, \spn{f_i}) &= (\tau
 \sigma, \spn{g_{\sigma(i)} \circ f_i})\text.
\end{align*}
\item \textbf{On vertical maps}: Given a functor $F \colon X \to Y$, we
give $SF \colon SX \to SY$ by
\begin{align*}
SF(n, \spn{x_i}) &= (n, \spn{Fx_i})\\
SF(\sigma, \spn{g_i}) &= (\sigma, \spn{Fg_i})\text.
\end{align*}
\item \textbf{On horizontal maps}: Given a profunctor $\X \colon X
  \tor Y$, we give the profunctor $S\X \colon SX \tor SY$ as
  follows:
\begin{itemize*}
\item \textbf{Proarrows} are
\[(\sigma, \spn{g_i}) \colon (n, \spn{y_i}) \tor (n, \spn{x_i})\text,\] where
$\sigma \in S_n$ and $g_i \colon y_i \tor x_{\sigma(i)}$ in $\X$ (no proarrows exist from $(n, \spn{y_i})$ to $(m, \langle{x_j}\rangle)$ for $n \neq m$);
\item \textbf{Actions} by maps of $SX$ and $SY$ are given in the obvious way, i.e.,
\[
(\tau, \spn{h_i}) \bullet (\sigma, \spn{g_i}) = (\tau
 \sigma, \spn{g_{\sigma(i)} \bullet f_i})
\]
for the left action by $SX$, and similarly for the right action by $SY$.
\end{itemize*}
\item \textbf{On cells}: Given a cell $\b F \colon \b X \Rightarrow \b Y$,
the cell $S \b F \colon S \b X \Rightarrow S \b Y$ is given by
\[S\b F\big((\sigma, \spn{g_i})\big) = (\sigma, \spn{\b F(g_i)})\text.\]
\end{itemize}
Vertical functoriality is immediate, whilst horizontal
pseudo-functoriality is easily defined and checked to be coherent.
There are straightforward variations on the above theme; we can
construct homomorphisms $T$ and $P \colon \dCat \to \dCat$ which lift,
respectively, the 2-functors on $\cat{Cat}$ for the `free
(non-symmetric) strict monoidal category' and the `free category with
finite products'.

There are general principles at work here: in all three cases, we have
a 2-functor $F$ on $\cat{Cat}$ which \emph{lifts} to a homomorphism
$\hat F$ on $\cat{Mod}$ in the sense of \cite{tanaka:thesis}. Any such
lifting will give rise to a double homomorphism on $\dCat$ which
`looks like' $F$ in the vertical direction and `looks like' $\hat F$
in the horizontal direction.
\end{Ex}

\begin{Defn}
Given morphisms $F, G \colon \db K \to \db L$ of pseudo double categories, a
\defn{vertical transformation} $\alpha \colon F \To G$ consists of
natural transformations $\alpha_0 \colon F_0 \Rightarrow G_0$ and $\alpha_1
\colon F_1 \Rightarrow G_1$ (and again, we shall use `$\alpha$' indifferently
for $\alpha_0$ and $\alpha_1$), subject to four straightforward axioms.

\end{Defn}

Given pseudo double categories $\db K$ and $\db L$, the double
morphisms $\db K \to \db L$ and vertical transformations between them
form a category $[\db K, \db L]_v$. These categories provide us with
hom-categories enriching $\cat{DblCat}$ to a 2-category.  Horizontal
composition of 2-cells is given by the horizontal composition in
$\cat{Cat}$ of the underlying natural transformations.

$[\db K, \db L]_v$ has a full subcategory $[\db K, \db L]_{v \psi}$
given by restricting to the double \emph{homomorphisms}.  Since double
homomorphisms are closed under horizontal composition, these fit
together to give the locally full sub-2-category $\cat{DblCat}_\psi$
of $\cat{DblCat}$, consisting of pseudo double categories, double
homomorphisms and vertical transformations.

\begin{Ex}\label{etamu}
  Following on from Examples \ref{dcat} and \ref{S}, we give a vertical transformation $\eta \colon \id_{\dCat}
  \Rightarrow S$ as follows. Its component at an object $X$ of
  $\dCat_0$ is the functor $\eta_X \colon X \rightarrow SX$ given by
\[\eta_X(x) = (1, \spn{x}) \quad \text{and} \quad \eta_X(f) = (\id_1, \spn{f})\text,\]
whilst its component at an object $\X$ of $\dCat_1$ is the cell $\eta_{\X} \colon \X \Rightarrow S\X$ given by
\[\eta_\X(g) = (\id_1, \spn{g})\text.\]
Likewise, we can give a vertical transformation $\mu \colon SS
\Rightarrow S$ which `flattens lists of lists' by removing the inner
sets of brackets.  It's easy to check that $\eta$ and $\mu$ as defined
above obey the monad laws
\[
  \mu \circ \eta S = \id_S = \mu \circ S \eta
    \quad \text{and} \quad
  \mu \circ \mu S  = \mu \circ S \mu
\]
and thus describe a monad on the object $\dCat$ in the 2-category
$\cat{DblCat}_\psi$, one which lifts the 2-monad for symmetric strict
monoidal categories on $\cat{Cat}$.

We can repeat the above exercise for the 2-monads for (non-symmetric)
strict monoidal categories and categories with finite products,
lifting them to double monads on $\dCat$.  Again, there are general
principles at work: we are utilising a \emph{pseudo-distributive law}
in the sense of \cite{marm:psd1, tanaka:thesis}, which allows us to
lift our 2-monad on $\cat{Cat}$ to a pseudomonad on $\cat{Mod}$. From
this, we can deduce the existence of a double monad on $\dCat$
combining the two.
\end{Ex}

In general, we shall call a monad in $\cat{DblCat}_\psi$ a \emph{double monad}:
Grandis and Par\'e consider such double monads and their more general
cousins, monads in $\cat{DblCat}$, in
\cite{grandispare:adjoints}.

\begin{Defn}
Given double morphisms $A_s, A_t \colon \db K \to \db L$, a
\defn{horizontal transformation} $\A \colon A_s \Tor A_t$ consists of
a \defn{components functor} $A_c \colon K_0 \to L_1$ (and to simplify
notation, we shall write $\A X$ for $A_cX$ and $\b A f$ for $A_c
f$) together with
special invertible maps
$A_\X \colon A_t\X \otimes \A{X_s} \to \A{X_t} \otimes A_s\X$ natural in $\X$,
which we call the \emph{pseudonaturality} of $\b A$; in pasting notation
\[\cd[@+0.5em]{
  A_sX_s \ar[r]|-{\object@{|}}^{A_s\X} \utwocell{dr}{A_\X} \ar[d]|-{\object@{|}}_{\A{X_s}} &
  A_sX_t \ar[d]|-{\object@{|}}^{\A{X_t}} \\
  A_tX_s \ar[r]|-{\object@{|}}_{A_t\X} &
  A_tX_t\text.
}\]
These data satisfy four evident axioms.
\end{Defn}

\begin{Ex}
The vertical transformations $\eta \colon \id_{\dCat} \Rightarrow S$
and $\mu \colon SS \Rightarrow S$ of Example \ref{etamu} have
horizontal counterparts $(\eta)_\ast \colon \id_{\dCat} \Tor S$ and
$(\mu)_\ast \colon SS \Tor S$, with components at $X \in (\dCat)_0$
given by
\[(\eta_\ast)_X = (\eta_X)_\ast \colon X \tor SX \quad \text{and} \quad
(\mu_\ast)_X = (\mu_X)_\ast \colon SSX \tor SX\text, \] where $(\
)_\ast$ is the usual embedding homomorphism $\cat{Cat} \to \cat{Mod}$.
We leave the remaining details to the reader.
\end{Ex}

\begin{Defn}
  Given horizontal transformations $\A \colon A_s \Tor A_t$ and $\b
  B \colon B_s \Tor B_t$, a \defn{modification} $\b \gamma \colon \A
  \Rrightarrow \B$ consists of
  a pair of vertical transformations
  $\gamma_s \colon A_s \To B_s$ (the `vertical source') and $\gamma_t \colon A_t \To B_t$ (the `vertical target');
  together with a natural transformation $\gamma_c \colon A_c \Rightarrow B_c$ (the `central natural
  transformation'). To simplify notation, we shall refer to the components of $\gamma_c$ as `the components of $\b
  \gamma$', and write a typical such component as $\b \gamma_X$. This data must
  satisfy three evident axioms.

We shall notate such a modification as:
\[\cd{
A_s \ar@{=>}[d]_{\gamma_s} \ar@{=>}[r]^{\A}|-{\object@{|}} \dthreecell{dr}{\b \gamma}& A_t \ar@{=>}[d]^{\gamma_t}  \\
B_s \ar@{=>}[r]|-{\object@{|}}_{\B} & B_t\text. }\]
\end{Defn}

Given two pseudo double categories $\db K$ and $\db L$, the horizontal
transformations and modifications between them form a category $[\db K, \db
L]_h$; further, there are two evident projections $s, t \colon [\db K, \db L]_h
\to [\db K, \db L]_v$ which provide data for a \defn{functor double category} $[\db K, \db L]$
as follows:
\begin{itemize}
\item The horizontal composite $(\C \colon C_s \Tor C_t) \otimes (\A \colon A_s \Tor
C_s)$ has components functor $C_c(\thg) \otimes A_c(\thg)$ and pseudonaturality
maps
\[(C \otimes A)_\X \colon C_t\X \otimes (\C{X_s} \otimes \A{X_s}) \to (\C{X_t} \otimes \A{X_t}) \otimes A_s\X\]
given by the pasting
\[\cd[@+0.5em]{
  A_sX_s \ar[r]|-{\object@{|}}^{A_s\X} \utwocell{dr}{A_\X} \ar[d]|-{\object@{|}}_{\A{X_s}} &
  A_sX_t \ar[d]|-{\object@{|}}^{\A{X_t}} \\
  C_sX_s \ar[r]|-{\object@{|}}^{C_s\X} \utwocell{dr}{C_\X} \ar[d]|-{\object@{|}}_{\C{X_s}} &
  C_sX_t \ar[d]|-{\object@{|}}^{\C{X_t}} \\
  C_tX_s \ar[r]|-{\object@{|}}_{C_t\X} &
  C_tX_t\text.
}\] Given modifications
\[\cd{
  A_s \ar@{=>}[d]_{\gamma_s} \ar@{=>}[r]^{\A}|-{\object@{|}} \dthreecell{dr}{\b \gamma}&
  C_s \ar@{=>}[d]^{\delta_s}\\
  B_s \ar@{=>}[r]|-{\object@{|}}_{\B} &
  D_s}
  \quad \text{and} \quad
  \cd{
  C_s \ar@{=>}[d]^{\delta_s} \ar@{=>}[r]^{\b C}|-{\object@{|}} \dthreecell{dr}{\b \delta}&
  C_t \ar@{=>}[d]^{\delta_t}  \\
  D_s \ar@{=>}[r]|-{\object@{|}}_{\b D} &
  D_t\text, }
\]
the composite modification $\b \delta \otimes \b \gamma$ has $(\delta \otimes
\gamma)_s = \gamma_s$, $(\delta \otimes \gamma)_t = \delta_t$ and component at
$X$ given by
$\b \delta_X \otimes \b \gamma_X \colon \b CX \otimes \b AX \to \b DX \otimes
\b BX$.
\item The horizontal unit $\I_F \colon F \Tor F$ at $F$ has components functor $\I_{F(\thg)}$, and pseudonaturality maps
$(I_F)_\X$ given by
\[(I_F)_\X = F\X \otimes \I_{FX_s} \xrightarrow{\r^{-1}_{F\X}} F\X \xrightarrow{\l_{F\X}} \I_{FX_t} \otimes F\X\text.\]
Given a vertical transformation $\alpha \colon F \To G$, the modification
$\I_\alpha$ has $(\I_\alpha)_s = \alpha = (\I_\alpha)_t$, and component at $X$
given by
$\I_{\alpha_X} \colon \I_{FX} \to \I_{GX}$.
\item Unit and associativity constraints $\l$, $\r$ and $\a$ for $[\db K, \db L]$ are given `componentwise' from those in $\db L$.
\end{itemize}
\noindent
 There is a sub-pseudo double category $[\db K, \db
L]_\psi$, given by restricting to \emph{homo}morphisms as objects,
and taking all vertical and horizontal transformations and
modifications between them.

\subsection{Whiskering of homomorphisms}\label{whisker}
Given a double morphism $G \colon \db L \to \db M$, we know by virtue of the
2-category structure of $\cat{DblCat}$ that we can `whisker' $F$ on either
side; that is, given vertical transformations
\[\alpha \colon F_1 \To F_2 \colon \db K \to \db L \quad \text{and} \quad
\beta \colon H_1 \To H_2 \colon \db M \to \db N\] we can form
vertical transformations
\[G\alpha \colon GF_1 \To GF_2 \colon \db K \to \db M \quad \text{and} \quad
\beta G \colon H_1G \To H_2G \colon \db L \to \db N\text.\] What
we shall do in this section is to produce a similar whiskering
operation on \emph{horizontal} transformations, and show that it
is compatible with the vertical whiskering:
\begin{Prop}
Let $G \colon \db L \to \db M$ be a double morphism. Then
`precomposition with $G$' extends to a strict double homomorphism
\[(\thg)G \colon [\db M, \db N] \to [\db L, \db N]\text.\]
\end{Prop}
\begin{proof}
We give $(\thg)G$ as follows:
\begin{itemize}
\item $\big((\thg)G\big)_0 \colon [\db M, \db N]_v \to [\db L, \db N]_v$ is given by
the whiskering operation in the 2-category $\cat{DblCat}$. Thus we
take the double morphism $H \colon \db M \to \db N$ to the double
morphism $HG \colon \db L \to \db N$ and the vertical
transformation $\alpha \colon H \To H'$ to the vertical
transformation $\alpha G \colon HG \To H'G$.

\item $\big((\thg)G\big)_1 \colon [\db M, \db N]_h \to [\db L, \db N]_h$ is given as
follows. Given a horizontal transformation $\A \colon A_s \Tor A_t$, the
horizontal transformation $\A G \colon A_sG \Tor A_tG$ has components functor
$A_cG_0$ (and therefore component at $X$ given by $\A G X \colon A_s GX \tor
A_t GX$) and pseudonaturality maps given by
\[(\A G)_{\X} = A_tG\X \otimes \A{GX_s} \xrightarrow{A_{G\X}} \A G{X_t} \otimes A_sG\X\text.\]
Given a modification $\b \gamma \colon \b A \Rrightarrow \b B$,

the
modification $\b \gamma G$ has $(\b \gamma G)_s = \gamma_s G$,
$(\b \gamma G)_t = \gamma_t G$, and

$(\b \gamma G)_c =\gamma_c G_0$, and therefore component at $X$
given by:
\[(\b \gamma G)_X = \b \gamma_{GX} \colon \A GX \to \B GX\text.\]
\end{itemize}
Visibly, $\big((\thg)G\big)_1$ and $\big((\thg)G\big)_0$ are compatible with
source and target, and we observe that $(\A \otimes \B)G = \A G \otimes \B G$
and $\I_HG = \I_{HG}$, so that $(\thg)G$ is a \emph{strict} homomorphism.
\end{proof}
\noindent We now move on to whiskerings on the left. As for bicategories, we
cannot in general whisker \emph{morphisms} with horizontal transformations on
the left; we must instead restrict to \emph{homo}morphisms.
\begin{Prop}
Let $G \colon \db L \to \db M$ be a double homomorphism. Then
`postcomposition with $G$' induces a double homomorphism
\[G(\thg) \colon [\db K, \db L] \to [\db K, \db M]\text.\]
\end{Prop}
\begin{proof}
We give $G(\thg)$ as follows:
\begin{itemize}
\item $\big(G(\thg)\big)_0 \colon [\db K, \db L]_v \to [\db K, \db M]_v$ is given by
the whiskering operation in the 2-category $\cat{DblCat}$. Thus we
take the double morphism $F \colon \db K \to \db L$ to the double
morphism $GF \colon \db K \to \db M$ and the vertical
transformation $\alpha \colon F \To F'$ to the vertical
transformation $G\alpha  \colon GF \To GF'$.
\item $\big(G(\thg)\big)_1 \colon [\db K, \db L]_h \to [\db K, \db M]_h$ is given as
follows. Given a horizontal transformation $\A \colon A_s \Tor A_t$, the
horizontal transformation $G \A \colon GA_s \Tor GA_t$ has components functor
$G_1 A_c$ (and therefore component at $X$ given by $G\A X \colon GA_sX \tor
GA_t X$) and pseudonaturality maps $(G\A)_\X$ given by

{\fontsize{8pt}{12pt}\selectfont \[GA_t\X \otimes G\A{X_s}
\xrightarrow{\m_{A_t\X, \A X_s}} G(A_t\X \otimes \A X_s)
\xrightarrow{GA_\X} G(\A{X_t} \otimes A_s\X)
\xrightarrow{\m^{-1}_{\A X_t, A_s\X}} G\A{X_t} \otimes
GA_s\X\text.\]}Given a modification $\b \gamma \colon \b A
\Rrightarrow \b B$, the modification $G \b \gamma$ has $(G\b
\gamma)_s = G\gamma_s$, $(G \b \gamma)_t = G \gamma_t$ and
$(G\b
\gamma)_c = G_1 \gamma_c$, and therefore component at $X$ given by
\[(G \b \gamma)_X = G\b \gamma_X \colon G\A X \to G\B X\text.\]
\end{itemize}
Again, it's clear that these functors are compatible with source and target. It
remains to give $\m$ and $\e$, so we take the special invertible modification
$\e_A \colon \I_{GA} \Rrightarrow G\I_A$ to have components
\[(\e_A)_X = \e_{AX} \colon \I_{GAX} \to G\I_{AX}\text.\]
and the special invertible modification
$\m_{\A, \B} \colon G\A \otimes G\B \Rrightarrow G(\A \otimes \B)$
to have components
\[(\m_{\A, \B})_X = \m_{\A X, \B X} \colon G\A X \otimes G\B X \to G(\A X \otimes
\B X)\text.\] Checking naturality and coherence is routine.
\end{proof}
\noindent Observe also that $G(\thg)$ and $(\thg)G$ restrict to respective
homomorphisms
\[(\thg)G \colon [\db M, \db N]_\psi \to [\db L, \db N]_\psi
\quad \text{ and } \quad G(\thg) \colon [\db K, \db L]_\psi \to [\db K, \db
M]_\psi\text.\] These propositions give us an `action' of homomorphisms on
functor pseudo double categories (we shall see below the precise sense in which
this \emph{is} an action), which can be extended from homomorphisms to
the vertical transformations between them. We begin with whiskerings on the right.
\begin{Prop}
Let $G$ and $G' \colon \db L \to \db M$ be double morphisms, and
let $\alpha \colon G \To G'$ be a vertical transformation. Then
precomposition with $\alpha$ induces a vertical transformation
\[
(\thg)\alpha \colon (\thg)G \Rightarrow (\thg)G' \colon [\db M,
\db N] \to [\db L, \db N]\text.\]
\end{Prop}
\begin{proof}
We give $(\thg)\alpha$ as follows:
\begin{itemize}
\item $\big((\thg)\alpha\big)_0$ has component at $H \in [\db M,
\db N]_v$ given by the map $H\alpha \colon HG \To HG'$ in $[\db L, \db N]_v$.
The naturality of these components in $H$ is the equality $\beta G' \circ
H\alpha = H'\alpha \circ \beta G$ in $\cat{DblCat}$;
\item $\big((\thg)\alpha\big)_1$ is given as follows. Its component at $\A \in [\db M, \db
N]_h$ is the modification $\b A \alpha \colon \b A G \Rrightarrow \b
A G'$ whose central natural transformation is $A_c\alpha_0$. The
naturality of these components in $\A$ follows from the equality
$\beta_c G'_0 \circ A_c\alpha_0 = A'_c\alpha_0 \circ \beta_c G_0$ in
$\cat{Cat}$.
\end{itemize}
These natural transformations are compatible with source and target, and checking
the vertical transformation axioms is routine.
\end{proof}
\begin{Prop}
Let $G$ and $G' \colon \db L \to \db M$ be double homomorphisms,
and let $\alpha \colon G \To G'$ be a vertical transformation.
Then postcomposition with $\alpha$ induces a vertical
transformation
\[
\alpha(\thg) \colon G(\thg) \Rightarrow G'(\thg) \colon [\db K,
\db L] \to [\db K, \db M]\text.
\]
\end{Prop}
\begin{proof}
We give the vertical transformation $\alpha(\thg)$ as follows:
\begin{itemize}
\item $\big(\alpha(\thg)\big)_0$ has component at $F \in [\db K,
\db L]_v$ given by the map $\alpha F \colon GF \To G'F$ in $[\db K, \db M]_v$.
The naturality of these components in $F$ is the equality  $G'\beta \circ
\alpha F = \alpha F' \circ G \beta$ in $\cat{DblCat}$.
\item $\big(\alpha(\thg)\big)_1$ has component at $\A \in [\db K, \db
L]_h$ given by the modification $\alpha \b \A \colon G \b A \Rrightarrow G' \b
A$ whose central natural transformation is $\alpha_1 A_c$. The naturality of
these components in $\A$ is the equality $G'_1\beta_c \circ \alpha_1 A_c  =
\alpha_1 A'_c \circ G_1 \beta_c$ in $\cat{Cat}$.
\end{itemize}
These natural transformations are compatible with source and target, and
checking coherence is routine.
\end{proof}
\noindent Observe that $\alpha(\thg)$ and $(\thg)\alpha$ restrict to respective
vertical transformations
\begin{align*}
(\thg)\alpha \colon (\thg)G & \Rightarrow (\thg)G' \colon [\db M, \db N]_\psi \to [\db L, \db N]_\psi\\
\alpha(\thg) \colon G(\thg) & \Rightarrow G'(\thg) \colon [\db K, \db L]_\psi \to [\db K, \db M]_\psi\text.
\end{align*}
We make one final remark: given a vertical transformation $\alpha
\colon G \To G'$ in $[\db L, \db M]_\psi$ and a modification $\b
\gamma \colon \b A \Rrightarrow \b B$ in $[\db M, \db N]$, the two
modifications $\b B \alpha \circ \b \gamma G$ and $\b \gamma G'
\circ \b A \alpha$ are the same, by naturality of
$\big((\thg)\alpha\big)_1$. Thus we shall write this common value as
$\b \gamma \alpha$. Similarly, if we have $\b \gamma \colon \b A
\Rrightarrow \b B$ now in $[\db K, \db L]$ we write $\alpha \b
\gamma$ for the modification $\alpha \B \circ G \b \gamma = G' \b
\gamma \circ \alpha \A$ in $[\db K, \db M]$.

\subsection{The hom 2-functor on $\cat{DblCat}_\psi$}
It's not hard to see that the operations of the previous section are functorial
with respect to vertical transformations. To be more precise, given double
categories $\db K$, $\db L$, $\db M$ and $\db N$, the above operations induce
functors
\begin{align*}
  [\db K, \thg] & \colon [\db L, \db M]_{v\psi} \to \big[\, [\db K, \db L], [\db K, \db M]\, \big]_{v \psi}
\\
\text{and }  [\thg, \db N] & \colon [\db L, \db M]_{v \psi} \to \big[\, [\db M, \db N], [\db L, \db N]\, \big]_{v
  \psi}\text,
\end{align*}
along with their `pseudo' restrictions
\begin{align*}
[\db K, \thg]_\psi & \colon [\db L, \db M]_{v\psi} \to \big[\, [\db K, \db L]_\psi, [\db K, \db M]_\psi\, \big]_{v \psi}\\
\text{and }  [\thg, \db N]_\psi & \colon [\db L, \db M]_{v \psi} \to \big[\,[\db M, \db N]_\psi, [\db L, \db N]_\psi\, \big]_{v  \psi}\text.
\end{align*}
Moreover, it's straightforward to check that the following
equalities hold:
\begin{align*}
 \big((\thg)G_1\big)G_2 &= (\thg)(G_1G_2)\text,
 &
 \big((\thg)\alpha_1\big)\alpha_2 &= (\thg)(\alpha_1\alpha_2)\text,
 \\
 G_1\big(G_2(\thg)\big) & = (G_1G_2)(\thg)\text, &
 \alpha_1\big(\alpha_2(\thg)\big) & = (\alpha_1\alpha_2)(\thg)\text, \\
 \big(G_1(\thg)\big)G_2 &= G_1\big((\thg)G_2\big)\text,&
 \text{and} \quad \big(\alpha_1 (\thg)\big) \alpha_2 &= \alpha_1 \big((\thg)\alpha_2\big)\text.
\end{align*}
which can be more succinctly stated as follows:
\begin{Prop}\label{internalhom} The functors $[\db K, \thg]$ and
$[\thg, \db N]$ defined above provide data for 2-functors
\[
[\db K, \thg] \colon \cat{DblCat}_\psi \to \cat{DblCat}_\psi\quad
\text{and } \quad[\thg, \db N] \colon \cat{DblCat}_\psi^\op \to
\cat{DblCat}_\psi
\]
which are compatible in the sense that they provide data for a
2-functor
\[[\,\thg, ?\,] \colon \cat{DblCat}_\psi^\op \times
\cat{DblCat}_\psi \to \cat{DblCat}_\psi\text.\] Similarly, the
functors $[\db K, \thg]_\psi$ and $[\thg, \db N]_\psi$ defined
above provide data for 2-functors
\[
[\db K, \thg]_\psi \colon \cat{DblCat}_\psi \to \cat{DblCat}_\psi
\quad \text{and} \quad [\thg, \db N]_\psi \colon \cat{DblCat}_\psi^\op \to
\cat{DblCat}_\psi
\]
which are compatible in the sense that they provide data for a
2-functor
\[[\,\thg, ?\,]_\psi \colon \cat{DblCat}_\psi^\op \times
\cat{DblCat}_\psi \to \cat{DblCat}_\psi\text.\]
\end{Prop}
\noindent Now, what \emph{are} these 2-functors? Does either of
the bivariant 2-functors provide an `internal hom' for
$\cat{DblCat}_\psi$? Let us make this question precise: observe
that $\cat{DblCat}_\psi$ has all finite products, and thus can be
viewed as a monoidal bicategory, with the tensor product given by
cartesian product. Then by an `internal hom' for
$\cat{DblCat}_\psi$, we mean a homomorphism of bicategories
\[\spn{\,\thg, ?\,} \colon \cat{DblCat}_\psi^\op \times
\cat{DblCat}_\psi \to \cat{DblCat}_\psi\] such that for all pseudo
double categories $\db K$, we have a biadjunction
$(\thg) \times \db K \dashv \spn{\db K, \thg}$.
In other words, $\spn{\,\thg, ?\,}$, if it exists, exhibits
$\cat{DblCat}_\psi$ as a \emph{biclosed} monoidal bicategory in
the sense of \cite{daystreet:monoidalbicats}.

Now, there is \emph{no} good biadjunction for the `lax hom'
2-functor $[\,\thg, ?\,]$, for the same reason as there is no good
whiskering on the left by morphisms: at some point, we have to produce
pseudo-naturality data for a horizontal transformation, and, due to
the laxity of the morphisms involved, no choice of such data exists.
However, it \emph{is} the case that the `pseudo hom' 2-functor
$[\,\thg, ?\,]_\psi$ provides an internal hom in the above
described sense. We don't intend to work through the rather messy
details here, but we do note that although both $(\thg) \times \db K$
and $[\,\db K, \thg\,]$ are 2-functors, the adjunction between them is
still only a \emph{bi}adjunction rather than an honest 2-adjunction.

%% file: chapter2.tex
\section{Clubs}\label{chap:cart1}
We now recall some of the basic definitions and results of the
theory of clubs.  A rather more detailed account of this material
can be found in \cite{kelly:clubs} or \cite{weber:clubs}.

\begin{Defn}\label{cartnat}
A natural transformation $\alpha \colon A \To S \colon \cat{C} \to \cat{D}$ is
called a \defn{cartesian natural transformation} if all its naturality squares
are pullbacks.
\end{Defn}
\noindent The following is immediate by elementary properties of pullback:
\begin{Prop}
Suppose that $\cat{C}$ has a terminal object $1$. Then a natural
transformation $\alpha \colon A \To S \colon \cat{C} \to \cat{D}$
is cartesian if and only if every naturality square of the form
\[\cd{
  AX
    \ar[r]^{A!}
    \ar[d]_{\alpha_X} &
  A1
    \ar[d]^{\alpha_1} \\
  SX
    \ar[r]_{S!} &
  S1
}\] is a pullback.
\end{Prop}
\noindent Thus, if we are given $S$, the cartesian
natural transformations into it are determined up to isomorphism by
their component $\alpha_1 \colon A1 \to S1$. We can make this
statement precise as follows. Given a category $\cat C$ and an object
$X \in \cat C$, the slice category $\cat C / X$ has:
\begin{itemize*}
\item \textbf{Objects} being pairs $(U, f)$ where $U \in \cat C$ and $f \colon U \to X$;
\item \textbf{Maps} $j \colon (U, f) \to (V, g)$ being maps $j \colon U \to V$ in $\cat C$ with $gj = f$.
\end{itemize*}
In particular, given a functor $S \colon \cat{C} \to \cat{D}$, we form the
slice category $[\cat C, \cat D] / S$; consider now the full subcategory of
this given by the objects $(A, \alpha)$ where $\alpha$ is a \emph{cartesian}
natural transformations into $S$. We write $Coll(S)$ for this subcategory and
call it the \defn{category of collections} over $S$. We have a functor $F
\colon Coll(S) \to \cat{D} / S1$ which evaluates at $1$:
\begin{align*}
F \colon Coll(S) & \to \cat{D} / S1\\
(A, \alpha) & \mapsto (A1, \alpha_1)\\
\gamma & \mapsto \gamma_1\text,
\end{align*}
and our above statement now becomes:
\begin{Prop}\label{cartplain}\cite{kelly:clubs}
Suppose $\cat{D}$ has all pullbacks; then evaluation at $1$ induces
an equivalence of categories $Coll(S) \simeq \cat{D} / S1$.
\end{Prop}
\noindent Now suppose we are given a category $\cat{C}$ together
with a monad $(S, \eta, \mu)$ on $\cat{C}$. As above, we can form
the slice category $[\cat C, \cat C] / S$, but now we can go
further; indeed, $[\cat C, \cat C]$ is a (strict) monoidal category
and $(S, \eta, \mu)$ is a monoid in it. Thus the slice category
$[\cat C, \cat C] / S$ acquires a canonical monoidal structure,
given by
\[\I = (\cdl{\id_{\cat C} \ar@2[r]^-{\eta} & S}) \quad \text{and} \quad
(A, \alpha) \otimes (B, \beta) = (\cdl{A B \ar@2[r]^{\alpha
    \beta} & SS \ar@2[r]^{\mu} & S})\text.\] This
structure is `canonical' in the following sense: giving a monoid $S$ in $[\cat
C, \cat C]$ is equivalent to giving a lax monoidal functor $\elt S \colon 1 \to
[\cat C, \cat C]$, and $[\cat C, \cat C] / S$ equipped with the above monoidal
structure is a \emph{lax limit} for this arrow in the 2-category of monoidal
categories, lax monoidal functors and lax monoidal transformations.

Now, we may naturally ask whether the subcategory $Coll(S)$ of $[\cat C, \cat
C] / S$ is closed under the above monoidal structure. Explicitly:
\begin{Defn}
We say that a subcategory $\cat D$ of a monoidal category $\cat C$ is a
\defn{monoidal subcategory} if $\cat D$ can be made into a
monoidal category such that the inclusion $\cat D \hookrightarrow
\cat C$ is a strict monoidal functor.
\end{Defn}
\begin{Defn}\label{clubdef}
We say that a monad $(S, \eta, \mu)$ is a \defn{club} on $\cat{C}$
if $Coll(S)$ is a monoidal subcategory of $[\cat C, \cat C] / S$.
\end{Defn}

\noindent Given a club $(S, \eta, \mu)$, we can exploit the equivalence of
categories $Coll(S) \simeq \cat C / S1$ to transport the monoidal
structure on $Coll(S)$ to a monoidal structure on $\cat C / S1$.
Explicitly, this monoidal structure has unit given by $\I = \eta_1 \colon 1
\to S1$, and tensor product $(a, \theta) \otimes (b, \phi)$ given by
the left-hand composite in the following diagram:
 \[\cd{
  a \otimes b
   \ar[r] \ar[d] &
  a
   \ar[d]^{\theta} \\
  Sb
   \ar[d]_{S\phi}
   \ar[r]_-{S!} &
  S1 \\
  SS1
   \ar[d]_{\mu_1} \\
  S1\text.
 }
 \]
Now, the above definition of club is not easy to work with in practice, so
the following alternative description is often useful:
\begin{Prop}\label{clubalt}\cite{kelly:clubs}
A monad $(S, \eta, \mu)$ is a club on $\cat C$ if and only if:
\begin{enumerate*}
\item $\eta$ is a cartesian natural transformation;
\item $\mu$ is a cartesian natural transformation;
\item $S$ preserves cartesian natural transformations into $S$:
that is, whenever $\alpha \colon A \Rightarrow S$ is cartesian, so is $S \alpha
\colon SA \Rightarrow SS$.
\end{enumerate*}
\end{Prop}
\begin{Ex}\label{clubsclubs}
Straightforward calculation using the previous proposition shows
all of the following to be clubs on $\cat{Cat}$:
\begin{itemize*}
\item The `free symmetric strict monoidal category' monad $S$;
\item The `free (non-symmetric) strict monoidal category' monad $T$;
\item The `free category with finite products' monad $P$.
\end{itemize*}
\end{Ex}
In Example \ref{etamu} of the previous section, we saw that $S$, $T$
and $P$ extend from 2-monads on $\cat{Cat}$ to double monads on
$\dCat$. What we are going to show is that $S$, $T$ and $P$ also
extend from \emph{clubs} on $\cat{Cat}$ to \emph{double clubs} on
$\dCat$.  To do this, we first need to know what we \emph{mean} by a
double club, and this is the objective of the next three sections.

%% file: chapter3.tex
\section{Double clubs I}

We shall assume without further mention that $\db K$
and $\db L$ are pseudo double categories such that:
\begin{itemize*}
\item $\db K$ has a
\defn{double terminal object}; that is, an object $1 \in K_0$
such that $1$ is terminal in $K_0$ and $\I_1$ is terminal in
$K_1$;
\item $L_1$ and $L_0$ have all pullbacks and are
equipped with a choice of such; and furthermore, $s$ and $t$ preserve these
choices \emph{strictly}.
\end{itemize*}
In the terminology of \cite{grandispare:limits}, this latter condition
amounts to a lax functorial choice of double pullbacks. In fact, we
can rephrase much of the work of this section \emph{globally}, in
terms of double pullbacks in double functor categories. However, by
doing so we would lose sight of why we have to impose technical
conditions such as property (hps) below. Therefore we shall work at
the \emph{local} level of components and leave it to the reader to
translate into a global view.

\begin{Ex}
  The pseudo double category $\dCat$ satisfies both the above criteria.
  The terminal category $1$ provides a double terminal object. For
  the lax functorial choice of double pullbacks, we observe that $\dCat_0 =
  \cat{Cat}$ certainly has all pullbacks, whilst
  $\dCat_1$ is isomorphic to the category $\cat{Cat} / \b 2$
  (where $\b 2$ is the arrow category $0 \rightarrow 1$), and hence
  also has all pullbacks. Further, given a choice of pullbacks in
  $\dCat_0$, we can choose pullbacks in $\dCat_1$ such that $s$ and
  $t$ strictly preserve them.
\end{Ex}

\subsection{Slice double categories}
\noindent We begin by extending the notion of slice category from plain
categories to double categories. The details of this construction are already
known, and can be found (along with a discussion of the more general `comma
double categories') in \cite{grandispare:adjoints}. Thus we shall merely recap
the details.

\begin{Defn}
A \defn{monad} in the pseudo double category $\db K$ consists of:
\begin{itemize*}
\item An object $X$ in $K_0$;
\item An object $\X \colon X \tor X$ in $K_1$;
\item Special maps
$\m \colon \X \otimes \X \to \X$ and $\e \colon \I_X \to \X$
subject to the commutativity of the usual unitality and associativity diagrams.
\end{itemize*}
\end{Defn}
\noindent
Equivalently, this is to give a double morphism $\elt {\b X} \colon 1
\to \db K$. So, given a pseudo double category $\db K$
together with a monad $(\b X, \m, \e)$ in $\db K$, we form the
\defn{slice double category} $\db K / \b X$ as follows: $(\db K /
\X)_1 = K_1 / \X$ and $(\db K / \X)_0 = K_0 / X$, whilst $s$ and $t$
are given by
\begin{align*}
 s(\b U \xrightarrow{\b f} \b X) &= (U_s \xrightarrow{f_s} X)\text, & s(\b j) &= j_s\text,\\
 t(\b U \xrightarrow{\b f} \b X) &= (U_t \xrightarrow{f_t} X) & \text{and } t(\b j) &= j_t\text.
\end{align*}
$\I$ and $\otimes$ are given on objects by
\[\I_{(U, f)} = (\I_U \xrightarrow{\I_f} \I_{X} \xrightarrow{\e} \X)
\quad \text{and} \quad (\b U, \b f) \otimes (\b V, \b g) = (\b U
\otimes \b V \xrightarrow{\b f \otimes \b g} \X \otimes \X
\xrightarrow{\m} \X)
\]
and inherit their action on maps from $\db K$, whilst
the natural transformations $\l$, $\r$ and $\a$  have components
inherited from $\db K$; that is,
\[
\l_{(\b U,  \b f)} = \l_{\b U}\text, \quad
 \r_{(\b U, \b f)} = \r_{\b U}\quad \text{and} \quad
 \a_{(\b U, \b f), (\b V, \b g), (\b W, \b h)} = \a_{\b U, \b V, \b W}\text.
\]
The remaining details are easily checked.  We now describe the slice
double categories we shall need for the theory of double clubs.
\begin{Prop}
Given a pseudo double category $\db K$ and an object $X \in K_0$,
the functor $\elt X \colon 1 \to K_0$ extends to a double
homomorphism $\elt{\I_X} \colon 1 \to \db K$.
\end{Prop}
\begin{proof}
To give $\elt{\I_X}$ is to give an `iso-monad' in $\db K$ whose
multiplication and unit are invertible; for this we take $\I_X
\colon X \tor X$, with multiplication and unit given by
\[\m = \l^{-1}_{\I_X} = \r^{-1}_{\I_X} \colon \I_X \otimes \I_X \to \I_X \quad \text{and} \quad \e = \id_{\I_X} \colon \I_X \to \I_X\text.\qedhere\]
\end{proof}
\noindent
In particular, given a double homomorphism $S \colon \db K \to \db L$, we have
the object $\id_{\db K} \in [\db K, \db K]_\psi$, and thus the double
homomorphism
\[1 \xrightarrow{\elt{\I_{\id_{\db K}}}} [\db K, \db K]_{\psi}
\xrightarrow{S(\thg)} [\db K, \db L]_{\psi}\text.\] Writing $S\I$ for the
corresponding monad in $[\db K, \db L]_\psi$, we can form the slice double category
$[\db K, \db L]_{\psi} / S\I$. Similarly, we have the monad $S\I_1$ given by
\[1 \xrightarrow{\elt{\I_1}} \db K
\xrightarrow{S} \db L\text,\] and so can form the slice double category $\db L
/ S\I_1$.

\begin{Ex}\label{catsi}
  Consider once more the double homomorphism $S \colon \dCat \to
  \dCat$ of Example \ref{S}. For this, the pseudo double category $\dCat / S\I_1$ has:
\begin{itemize}
\item \textbf{Objects} $(X, F)$ given by a category $X$ together
with a functor $F \colon X \to S1$. We observe that we can identify
$S1$ with (a skeleton of) the category of \emph{finite sets and
bijections}.
\item \textbf{Vertical maps} $H \colon (X, F) \to (Y, G)$ given by commutative
triangles
\[\cd{
  X
   \ar[rr]^{H}
   \ar[dr]_{F} & &
  Y
   \ar[dl]^{G} \\ &
  S1\text.
}
\]
\item \textbf{Horizontal maps} $(\b X, \b F) \colon (X_s, F_s) \tor (X_t, F_t)$ given by a profunctor $\X \colon X_s \tor X_t$ together with a cell
\[\cd{
  X_s \ar[d]_{F_s} \ar[r]^{\X}|-{\object@{|}} \dtwocell{dr}{\b F}& X_t \ar[d]^{F_t}  \\
  S1 \ar[r]|-{\object@{|}}_{S\I_1} & S1\text.}\] We identify the
profunctor $S\I_1 \colon S1 \tor S1$ with the hom functor on $S1$;
thus to give a horizontal map $(\b X, \b F)$ is to give a profunctor
$\X$ together with an assignation to each proarrow $f$ of $\X$ an
arrow $\b Ff$ of $S1$, compatible with $F_s$ and $F_t$.
\item \textbf{Cells} $\b H \colon (\b X, \b F) \Rightarrow (\b Y, \b G)$
are given by commutative triangles of cells in $\dCat$
\[\cd{
  \b X
   \ar@2[rr]^{\b H}
   \ar@2[dr]_{\b F} & &
  \b Y
   \ar@2[dl]^{\b G} \\ &
  S\I_1\text;
}
\]
thus to each proarrow $f$ of $\X$, we assign a compatible proarrow
$\b Hf$ of $\Y$ such that $\b G\b H f = \b Ff$.
\item \textbf{Horizontal identity} is given on objects $(X, F)$ by
  $(\I_X, \hat{\I}_F)$, where $\I_X$ is the identity profunctor on $X$
  and $\hat{\I}_F$ is given by $\hat{\I}_F(\I_f) = Ff$, for $f$ an
  arrow of $X$.
\item \textbf{Horizontal composition} is given by
$(\b X, \b F) \otimes (\b X', \b F') = (\b X \otimes \b X', \b F \hat{\otimes} \b F')$, where $\b X \otimes \b X'$ is usual profunctor
composition, and where $(\b F \hat{\otimes} \b F')(f \otimes f') = \b F f \circ \b F'(f')$.
\end{itemize}
\noindent The pseudo double category $[\dCat, \dCat]_\psi / S\I$ has:
\begin{itemize}
\item \textbf{Objects} $(A, \alpha)$ given by a double homomorphism $A \colon \dCat \to \dCat$
together with a vertical transformation $\alpha \colon A \To S$.
\item \textbf{Vertical maps} $\gamma \colon (A, \alpha) \to (B, \beta)$ given by commutative triangles
\[\cd{
  A
   \ar@2[rr]^{\gamma}
   \ar@2[dr]_{\alpha} & &
  B
   \ar@2[dl]^{\beta} \\ &
  S\text.
}
\]
\item \textbf{Horizontal maps} $(\b A, \b \alpha) \colon (A_s, \alpha_s) \tor (A_t, \alpha_t)$ given by pairs $(\b A, \b \alpha)$
where $\b A$ is a horizontal transformation and $\b \alpha$ a modification as
follows:
\[\cd{
A_s \ar@{=>}[d]_{\alpha_s} \ar@{=>}[r]^{\b A}|-{\object@{|}} \dthreecell{dr}{\b \alpha}& A_t \ar@{=>}[d]^{\alpha_t}  \\
S \ar@{=>}[r]|-{\object@{|}}_{S\I} & S\text;}\]
\item \textbf{Cells} $\b \gamma \colon (\b A, \b \alpha) \Rightarrow (\b B, \b \beta)$ given by
commutative triangles
\[\cd{
  \b A
   \ar@3[rr]^{\b \gamma}
   \ar@3[dr]_{\b \alpha} & &
  \b B
   \ar@3[dl]^{\b \beta} \\ &
  S\I\text.
}
\]
\item \textbf{Horizontal identities} given on objects $(A, \alpha)$ by
\[\I_{(A, \alpha)} = \cdl[@!]{\I_A \ar@3[r]^-{\I_{\alpha}} & \I_{S} \ar@3[r]^-{\e} & S\I}\]
(where $\e$ is the unit of the monad $S\I$, with components $\e_X \colon
\I_{SX} \to S\I_X$), and on maps $\gamma \colon (A, \alpha) \to (B, \beta)$ by
\[\cd{ \I_A
   \ar@3[rr]^{\I_\gamma}
   \ar@3[dr]_{\e \circ \I_\alpha} & &
  \b B
   \ar@3[dl]^{\e \circ \I_\beta} \\ &
  S\I\text;
}\]
\item \textbf{Horizontal composition} given on objects by
\[(\b A, \b \alpha) \otimes (\b A', \b \alpha') = \big( \cdl[@!@+1em]{\b A \otimes \b A' \ar@3[r]^-{\b \alpha \otimes \b \alpha'} & S\I \otimes S\I \ar@3[r]^-{\m} & S\I} \big)\]
(where $\m$ is the multiplication of the monad $S\I$, with components
\[\m_X = S\I_X \otimes S\I_X \xrightarrow{\m_{\I_X, \I_X}} S(\I_X \otimes \I_X) \xrightarrow{S\l^{-1}_{\I_X}} S\I_X\text{ ),}\]
and on maps by
\[\cd{ \b A \otimes \b A'
   \ar@3[rr]^{\b \gamma \otimes \b \gamma'}
   \ar@3[dr]_{\m \circ (\b \alpha \otimes \b \alpha')} & &
  \b B \otimes \b B'
   \ar@3[dl]^{\m \circ (\b \beta \otimes \b \beta')} \\ &
  S\I\text.
}\]
\end{itemize}
\end{Ex}

\subsection{The double category of collections}\label{cartrequire}
We return now to our general theory. We should like to restrict from
the full double slice category $[\db K, \db L]_\psi/S\I$ to
something mimicking the category of collections. To do this, we need
a double category analogue of Definition \ref{cartnat}'s `cartesian
natural transformation':
\begin{Defn}\hfill
\begin{itemize}
\item A vertical transformation $\alpha \colon F \To G \colon \db K \to \db L$ is
called a \defn{cartesian vertical transformation} if the natural
transformations $\alpha_1 \colon F_1 \Rightarrow G_1$ and
$\alpha_0 \colon F_0 \Rightarrow G_0$ are cartesian;
\item A modification $\b \gamma \colon \b A \Rrightarrow \b B$
is called a \defn{cartesian modification} if $\gamma_s$ and
$\gamma_t$ are cartesian vertical transformations and the natural
transformation $\gamma_c \colon A_c \Rightarrow B_c$ is cartesian.
\end{itemize}
\end{Defn}
\noindent We should like the double category of collections $\coll(S)$ to have:
\begin{itemize}
\item $\coll(S)_0$ being the full subcategory of $\big([\db K, \db
  L]_\psi/S\I\big)_0$ whose objects are the cartesian vertical
  transformations into $S$;
\item $\coll(S)_1$ being the full subcategory of $\big([\db K, \db L]_\psi/S\I\big)_1$
whose objects are the cartesian modifications into $S\I$,
\end{itemize}
with the remaining data inherited from the double category $[\db
K, \db L]_\psi/S\I$. In order for this to make sense, we need
$\coll(S)$ to be closed under the horizontal units and composition
of $[\db K, \db L]_\psi/S\I$, for which we require $S$ to have the
following property.
\begin{Defn}
Let $S \colon \db K \to \db L$ be a double homomorphism; we say
that $S$ has
\defn{property (hps)} (horizontal pullback stability) if it satisfies:
\begin{itemize}
\item \textbf{Property (hps1)}: given horizontally composable
pullbacks
\[
\cd{
{\b A} \ar[r]^{\b p_1} \ar[d]_{\b p_2} & {\b B} \ar[d]^{\b f} \\
{S\b C} \ar[r]_{S!} & {S \I_1} } \quad \text{and} \quad
\cd{
{\b A'} \ar[r]^{\b p'_1} \ar[d]_{\b p'_2} & {\b B'} \ar[d]^{\b f'} \\
{S\b C'} \ar[r]_{S!} & {S \I_1}\text, }
\]
in $L_1$, the diagram
\[
\cd{
{\b A' \otimes \b A} \ar[r]^{\b p'_1 \otimes \b p_1} \ar[d]_{\b p'_2 \otimes \b p_2} & {\b B' \otimes \b B} \ar[d]^{\b f' \otimes \b f} \\
{S \b C' \otimes S\b C} \ar[r]_{S! \otimes S!} & {S \I_1 \otimes S
\I_1} }
\] is a pullback in $L_1$; and
\item \textbf{Property (hps2)}: given a pullback
\[
\cd{
{A} \ar[r]^{p_1} \ar[d]_{p_2} & {B} \ar[d]^{f} \\
{SC} \ar[r]_{S!} & {S1} }
\]
in $L_0$, the diagram
\[
\cd{
{\I_A} \ar[r]^{\I_{p_1}} \ar[d]_{\I_{p_2}} & {\I_B} \ar[d]^{\I_f} \\
{\I_{SC}} \ar[r]_{\I_{S!}} & {\I_{S1}} }
\]
is a pullback in $L_1$.
\end{itemize}
\end{Defn}

\begin{Prop}
Given a homomorphism $S \colon \db K \to \db L$ with property
(hps), the categories $\coll(S)_0$ and $\coll(S)_1$ provide data
for a pseudo double category whose remaining data is inherited
from $[\db K, \db L]_\psi / S\I$.
\end{Prop}
\begin{proof}
We must check that the horizontal units of $[\db K, \db L]_\psi /
S\I$ are cartesian modifications, and that the horizontal
composition of two cartesian modifications is another cartesian
modification. For the first of these, given $(A, \alpha) \in
\coll(S)_0$, we have $\I_{(A, \alpha)}$ given by the modification
\[\I_{(A, \alpha)} = \cdl{\I_A \ar@3[r]^{\I_\alpha} & \I_S \ar@3[r]^{\e} & S\I}\text;\]
so consider the diagram
\[\cd{
{\I_{AX}} \ar[r]^{\I_{A!}} \ar[d]_{\I_{\alpha_X}} &
{\I_{A1}} \ar[d]^{\I_{\alpha_1}} \\
{\I_{SX}} \ar[r]^{\I_{S!}} \ar[d]_{\e_X} &
{\I_{S1}} \ar[d]^{\e_1} \\
{S\I_X} \ar[r]_{S\I_!} & {S\I_1}\text. }\] It follows from
property (hps2) and the cartesianness of $\alpha$ that the top
square is a pullback; and the lower square commutes, and so is a
pullback since both vertical arrows are isomorphisms. Thus the
outer edge is again a pullback, and so $\I_{(A, \alpha)}$ is cartesian as
required.

For the second, suppose we are given horizontally composable
objects $(\b A, \b \alpha)$ and $(\b B, \b \beta)$ of
$\coll(S)_1$; we must show that the modification
\[\cdl[@C+1em]{\b A \otimes \b B \ar@3{->}[r]^-{\b \alpha \otimes \b \beta} & S\I \otimes S\I \ar@3[r]^-{\m} & S\I}\]
is also cartesian. So consider the diagram:
\[
\cd[@+1em]{
\A X \otimes \B X \ar[r]^{\A ! \otimes \B !} \ar[d]_{\b \alpha_X \otimes \b \beta_X}& \A 1 \otimes \B 1 \ar[d]^{\b \alpha_1 \otimes \b \beta_1} \\
S\I_X \otimes S\I_X \ar[r]_{S\I_! \otimes S\I_!} \ar[d]_{\m_X} &
S\I_1 \otimes S\I_1 \ar[d]^{\m_1} \\
S\I_X \ar[r]_{S\I_!} & S\I_1\text. }
\]
The upper square is a pullback by property (hps2) and the
cartesianness of $\alpha$ and $\beta$; the lower square commutes
and has isomorphisms down the sides, and hence is a pullback. So
the outer edge is also a pullback as required.
\end{proof}

\subsection{Evaluation at $1$ in $\coll(S)$}
In order to see that our definition of $\coll(S)$ is the correct
one, we need to show that there is a suitable analogue at the pseudo
double category level of the equivalence of categories $Coll(S)
\simeq \cat D / S1$ exhibited in Proposition \ref{cartplain}.

For this, we need a suitable notion of `equivalence of double
categories'. There is an obvious candidate for this, namely
equivalence in the 2-category $\cat{DblCat}_\psi$, and the following
proposition gives us an elementary characterisation of such
equivalences.
\begin{Prop}\label{equiv}
Suppose we are given double categories $\db K$ and $\db L$, and:
\begin{itemize*}
\item A double homomorphism $F \colon \db K \to \db L$;
\item Functors $G_1 \colon L_1 \to K_1$ and $G_0 \colon L_0 \to K_0$;
\item Natural isomorphisms $\eta_i \colon \id_{K_i} \cong G_iF_i$ and
$\epsilon_i \colon F_iG_i \cong \id_{K_i}$ ($i = 0, 1$);
\end{itemize*}
such that $s G_1 = G_0 s$, $t G_1 = G_0 t$, $s \epsilon_1 = \epsilon_0
s$, $t \epsilon_1 = \epsilon_0 t$, $s \eta_1 = \eta_0 s$ and $t \eta_1
= \eta_0 t$. Then $\db K$ and $\db L$ are equivalent in
$\cat{DblCat}_\psi$.
\end{Prop}
\begin{proof}
See Appendix A, Corollary \ref{equivv}.
\end{proof}
\noindent Now let $S \colon \db K \to \db L$
be a double homomorphism with property (hps), and consider the double
category of collections $\coll(S)$. We have a strict homomorphism $F \colon
\coll(S) \to \db L / S\I_1$ which `evaluates at $1$':
\[
\begin{aligned}
F_0 \colon \coll(S)_0 & \to L_0 / S1\\
(A, \alpha) & \mapsto (A1, \alpha_1)\\
\gamma & \mapsto \gamma_1
\end{aligned}
\qquad \text{and} \qquad
\begin{aligned}
F_1 \colon \coll(S)_1 & \to L_1 / S\I_1\\
(\b A, \b \alpha) & \mapsto (\b A1, \b \alpha_1)\\
\b \gamma & \mapsto \b \gamma_1\text.
\end{aligned}
\]
\noindent Using this, we can prove the following analogue of Proposition \ref{cartplain}.
\begin{Prop}
  Let $S$ be a homomorphism $\db K \to \db L$ satisfying property
  (hps). Then evaluation at $1$ induces an equivalence of double
  categories $\coll(S) \simeq \ks$.
\end{Prop}
\begin{proof}
  We exhibit all the data required for Proposition \ref{equiv}. We
  have the strict homomorphism $F \colon \coll(S) \to \ks$ as above;
in the opposite direction, we must exhibit functors $G_i \colon
(\ks)_i \to \coll(S)_i$. We can form categories of collections
$Coll(S_0)$ and $Coll(S_1)$, and by Proposition \ref{cartplain} we
have equivalences of categories
\[Coll(S_0) \simeq L_0 / S1 \quad \text{and} \quad Coll(S_1) \simeq L_1 / S\I_1\]
where the rightward direction of these equivalences is given by
evaluation at $1$ and $\I_1$ respectively. We are now ready to
give $G_0$:
\begin{itemize}
\item \textbf{On objects}: given an object $(a, \theta) \in L_0 / S1$, under the first
equivalence we produce an object $(A_0, \alpha_0) \in Coll(S_0)$. We can also
form the object $\I_{(a, \theta)} \in L_1 / S\I_1$: under the second
equivalence this produces an object $(A_1, \alpha_1) \in Coll(S_1)$.
Explicitly, $A_0$, $\alpha_0$, $A_1$ and $\alpha_1$ are the specified objects
and maps in the following pullback diagrams:
\[
\cd{
 A_0X \ar[r]^-{A_0 !} \ar[d]_{(\alpha_0)_X} &
 a \ar[d]^{\theta} \\
 SX \ar[r]_-{S!} &
 S1
} \quad \text{ and } \quad
\cd{
 A_1\X \ar[r]^-{A_1 !} \ar[dd]_{(\alpha_1)_{\X}} &
 \I_a \ar[d]^{\I_\theta} \\ &
 \I_{S1} \ar[d]^{\e_1} \\
 S\X \ar[r]_-{S!} &
 S\I_1\text.
}
\]
Since $s$ and $t$ strictly preserve pullbacks, its easy to see that $A_1$ and
$A_0$, and similarly $\alpha_1$ and $\alpha_0$, are compatible with source and
target. We aim to equip $A = (A_0, A_1)$ with the structure of a double
homomorphism, and to show that $\alpha = (\alpha_0, \alpha_1)$ becomes a
cartesian vertical transformation with respect to this structure. To do this,
we must produce special natural isomorphisms
\[\m_{\X, \Y} \colon A\X \otimes A\Y \to A(\X \otimes \Y) \quad
\text{and} \quad \e_X \colon \I_{AX} \to A\I_X\text.\] So consider
the diagram:
\[
\cd{
 A\X \otimes A\Y
  \ar[rr]^{A! \otimes A!}
  \ar[dd]_{\alpha_\X \otimes \alpha_\Y} & &
 {\I_a \otimes \I_a}
  \ar[dr]^-{\r^{-1}_{\I_a}}
  \ar[dd]^(0.65){(\e_1 \circ \I_\theta) \otimes (\e_1 \circ \I_\theta)} \\ &
 A(\X \otimes \Y)
  \ar[dd]_(0.65){\alpha_{\X \otimes \Y}}
  \ar[rr]^{A!} & &
 {\I_a}
  \ar[dd]^{\e_1 \circ \I_\theta} \\
 S\X \otimes S\Y
  \ar[rr]^(0.65){S! \otimes S!}
  \ar[dr]_{\m_{\X, \Y}} & &
 S\I_1 \otimes S\I_1
  \ar[dr]^{\m_1} \\ &
 S(\X \otimes \Y)
  \ar[rr]_{S!} & &
 S\I_1
}
\]
The front face is a pullback by definition; the back face by
property (hps1). All the diagonal maps are isomorphisms, and the
bottom and right faces commute by the coherence axioms for $S$ and
$\db L$. Thus we induce a unique isomorphism $A\X \otimes A\Y \to
A(\X \otimes \Y)$ along the missing diagonal. Arguing identically
for the unit, we induce a unique isomorphism $\I_{AX} \to A\I_X$.
All required naturality and coherence now follows straightforwardly
using the existing coherence and the universal property of pullback.

\item \textbf{On maps}: suppose we have a map $\psi \colon (a, \theta) \to (b,
\phi)$ in $L_0 / S1$, with $G_0(a, \theta) = (A, \alpha)$ and
$G_0(b, \phi) = (B, \beta)$. Then we must produce a map $\gamma
\colon (A, \alpha) \to (B, \beta)$; that is, a vertical
transformation $\gamma \colon A \To B$ making the diagram
\[
\cd{ A \ar@2[rr]^{\gamma} \ar@2[dr]_{\alpha} & & B
\ar@2[dl]^{\beta} \\ & S}
\]
commute. Now, using the equivalences $L_0 / S1 \simeq Coll(S_0)$
and $L_1 / S\I_1 \simeq Coll(S_1)$ as before, we produce natural
transformations $\gamma_0$ and $\gamma_1$ making
\[
\cd{ A_0 \ar@2[rr]^{\gamma_0} \ar@2[dr]_{\alpha_0} & & B_0
\ar@2[dl]^{\beta_0} \\ & S_0} \quad \text{and} \quad \cd{ A_1
\ar@2[rr]^{\gamma_1} \ar@2[dr]_{\alpha_1} & & B_1
\ar@2[dl]^{\beta_1} \\ & S_1}
\]
commute. We aim to show that $\gamma = (\gamma_0, \gamma_1)$ becomes a vertical
transformation. Compatibility with source and target follows as before, whilst
the other two axioms follow from the naturality of $\r^{-1}$ and the universal
property of pullback.
\end{itemize}
We now move on to $G_1$. Suppose we have an object
\[\cd{
a_s \ar[d]_{\theta_s} \ar[r]^{\b a}|-{\object@{|}} \dtwocell{dr}{\b \theta}& a_t \ar[d]^{\theta_t}  \\
S1 \ar[r]|-{\object@{|}}_{S \I_1} & S1 }\] of $L_1 / S\I_1$, with
$G_0(a_s, \theta_s) = (A_s, \alpha_s)$ and $G_0(a_t, \theta_t) =
(A_t, \alpha_t)$, say. Then we must produce an object $(\b A, \b
\alpha) \in \coll(S)_1$ as follows:
\[\cd{
A_s \ar@2[d]_{\alpha_s} \ar@2[r]^{\b A}|-{\object@{|}} \dthreecell{dr}{\b \alpha}& A_t \ar@2[d]^{\alpha_t}  \\
S \ar@2[r]|-{\object@{|}}_{S\I} & S\text. }\] Under the equivalence $L_1 /
S\I_1 \simeq Coll(S_1)$, we take $(\b a, \b \theta)$ to a functor $A \colon K_1
\to L_1$ and a cartesian natural transformation $\alpha \colon A \Rightarrow
S_1$. Thus we specify the horizontal transformation $\b A$ to have source
$A_s$, target $A_t$ and components functor $A_c = A\I \colon K_0 \to L_1$.
Similarly, we take the modification $\b \alpha$ to have source $\alpha_s$,
target $\alpha_t$ and central natural transformation $\alpha_c = \alpha\I
\colon A\I \Rightarrow S\I \colon K_0 \to L_1$. Explicitly, $\A X$ and
$\alpha_X$ will be the indicated arrows in the following pullback diagram:
\[\cd{
\A X \ar[r]^{\A !} \ar[d]_{\alpha_X} & {\b a} \ar[d]^{\b \theta} \\
S\I_X \ar[r]_{S!} & S\I_1}\text.
\]
We must now specify the pseudonaturality maps for $\b A$. So consider the
diagram
\[
\cd{
 A_t\X \otimes \A X_s
  \ar[rr]^{A_t! \otimes \A !}
  \ar[dd]_{(\alpha_t)_{\X} \otimes {\b \alpha}_{X_s}} & &
 {\I_{a_t} \otimes \b a}
  \ar[dd]^(0.65){(\e_1 \circ \I_{\theta_t}) \otimes \b \theta}
  \ar[dr]^{\r_{\b a} \circ \l^{-1}_{\b a}} \\ &
 \A X_t \otimes A_s\X
  \ar[rr]^(0.35){\A ! \otimes A_s!}
  \ar[dd]_(0.35){{\b \alpha}_{X_t} \otimes (\alpha_s)_{\X}} & &
 {\b a \otimes \I_{a_s}}
  \ar[dd]^{\b \theta \otimes (\e_1 \circ \I_{\theta_s})} \\
 S\X \otimes S\I_{X_s}
  \ar[rr]^(0.65){S! \otimes S!}
  \ar[dr]_{(S\I)_{\X}} & &
 S\I_1 \otimes S\I_1
  \ar[dr]^{\id} \\ &
 S\I_{X_t} \otimes S\X
  \ar[rr]_{S! \otimes S!} & &
 S\I_1 \otimes S\I_1\text.
}\] The front and back faces are pullbacks by property (hps1) and
the diagonal maps are all isomorphisms. It's easy to check that the
bottom and right faces commute, and thus we induce a unique
isomorphism along the missing diagonal, which will be the
pseudonaturality map $A_\X$. Again, all required naturality and
coherence follows easily using the existing naturality and coherence
and the universal property of pullback.

We now give $G_1$ on maps. Given a map $\b \psi \colon (\b a, \b
\theta) \to (\b b, \b \phi)$ in $K_1 / S\I_1$, we must produce a
map $\b \gamma \colon (\A, \b \alpha) \to (\B, \b \beta)$ of
$\coll(S)_1$, and thus a modification $\b \gamma \colon \b A
\Rrightarrow \b B$ fitting into the diagram
\[\cd{\b A \ar@3[rr]^{\b \gamma} \ar@3[dr]_{\b \alpha} & & \b B \ar@3[dl]^{\b \beta} \\ & S\I\text.
}\] For its source and target, we take the vertical
transformations
\[
\gamma_s = G_0(\psi_s) \colon A_s \Rightarrow B_s \quad \text{and}
\quad \gamma_t = G_0(\psi_t) \colon A_t \Rightarrow B_t\text.\]
For the central natural transformation, we apply once more the
equivalence $L_1 / S\I_1 \cong Coll(S_1)$ to get a commuting
diagram
\[\cd{A \ar@2[rr]^{\gamma} \ar@2[dr]_{\alpha} & & B \ar@2[dl]^{\beta} \\ & S_1\text,
}\] in the functor category $[L_1, L_1]$. We need a natural transformation
$\gamma_c \colon A_c \Rightarrow B_c$, and from above we have $A_c = A\I$ and
$B_c = B\I$; so we take $\gamma_c = \gamma\I$. This this provides coherent data
for a modification follows by an argument similar to above. Finally, we note
that we have
\[\alpha_c = \alpha \I = (\beta \circ \gamma)\I = \beta \I \circ \gamma \I =
\beta_c \circ \gamma_c\] as required. This completes the
definition of $G_1$.

By construction, it is immediate that $t G_1 = G_0 t$ and $s G_1 = G_0
s$; so we need to show that $(F_0, G_0)$ and $(F_1, G_1)$ provide data
for equivalences of categories. First note that if we choose pullbacks
in $L_0$ and $L_1$ such that the pullback of identity arrows are
identity arrows then we have
\[F_0G_0 = \id_{L_0/S1} \quad \text{and} \quad F_1G_1 = \id_{L_1 / S\I_1}\text.\]
Conversely, it's an easy exercise using the universal property of
pullback to construct natural isomorphisms $\eta_0 \colon
\id_{\coll(S)_0} \Rightarrow G_0F_0$ and $\eta_1 \colon
\id_{\coll(S)_1} \Rightarrow G_1F_1$, and to show that they are
compatible with source and target maps as required. Thus we have all
the requirements for Corollary \ref{equiv}, and so have an
equivalence of double categories $\coll(S) \simeq \db K / S\I_1$.
\end{proof}


%% file: chapter4.tex
\section{Monoidal double categories}

To complete our exposition of the theory of double clubs, we need a
suitable generalisation of \emph{monoidal category} to the double
category level. This is fairly straightforward: recall that the
2-category $\cat{DblCat}_\psi$ has finite products, given in the
obvious way, and hence becomes a (cartesian) monoidal bicategory
\cite{gps}. Thus we can define
\begin{Defn}
A \defn{monoidal double category} is a pseudomonoid \cite{daystreet:monoidalbicats, mccrudden:thesis} in
$\cat{DblCat}_\psi$.
\end{Defn}
\noindent However, this definition is too abstract to work with in practice; we use instead the following alternative characterisation, the proof of which is entirely routine:
\begin{Prop}
  Giving a monoidal double category $\db K$ is equivalent to giving a
  double category $\db K$ such that
\begin{itemize*}
\item $K_0$ is a (not necessarily strict) monoidal category, with data $(\ten_0, \elt{\unit}, \alpha_0, \lambda_0, \rho_0)$;
\item $K_1$ is a (not necessarily strict) monoidal category, with data $(\ten_1, \b{\elt{\unit}}, \alpha_1, \lambda_1, \rho_1)$;
\item The functors $s$ and $t \colon K_1 \to K_0$ are strict
monoidal;
\item The functors $\b I \colon K_0 \to K_1$ and $\otimes \colon K_1 \mathbin{\vphantom{\times}_s\mathord{\times}_t} K_1 \to K_1$
are strong monoidal (where $K_1 \mathbin{\vphantom{\times}_s\mathord{\times}_t}
K_1$ acquires its monoidal structure via pullback along the strict monoidal
functors $s$ and $t$);
\item The associativity and unitality natural transformations
  $\a$, $\l$ and $\r$ for $\db K$ are monoidal natural
  transformations.
\end{itemize*}
\end{Prop}
\noindent We note that the data making $\ten$ and $e$ strong
monoidal amounts to giving invertible special maps in $K_1$ as follows:
\begin{gather*}
 k_{\W, \X, \Y, \Z} \colon (\W \bten \X) \otimes (\Y \bten \Z) \to (\W \otimes \Y) \bten (\X \otimes \Z)\text,\\
  u_{X, Y} \colon \I_{W \ten X} \to \I_W \bten \I_X\text,\\
  k_{\b{\unit}} \colon \b\unit \otimes \b\unit \to \b\unit \quad \text{and} \quad u_{\unit} \colon \I_{\unit} \to \b\unit\text,
\end{gather*}
natural in all variables and obeying a number of coherence diagrams.
\begin{Ex}
The pseudo double category $\dCat$ of Example \ref{dcat} becomes a
monoidal double category where $\ten$ is given on objects by
cartesian product of categories, extended in the evident way to
vertical maps, horizontal maps and cells. More generally, the pseudo
double category $\mathcal V$-$\dCat$ becomes a monoidal double
category where $\ten$ is now given by tensor product of $\mathcal
V$-categories.
\end{Ex}

We turn now to the apposite notion of \emph{map} between two monoidal
double categories. The obvious candidate is that of a \emph{lax map of
  pseudomonoids} \cite{daystreet:monoidalbicats} in
$\cat{DblCat}_\psi$. However, the underlying double morphism of such a
map is necessarily a \emph{homo}morphism, and this is not sufficiently
general.

To overcome this, we observe that the 2-category $\cat{DblCat}$
also has finite products, and that the inclusion $\cat{DblCat}_{\psi}
\to \cat{DblCat}$ preserves them. So we view a monoidal double
category \emph{a fortiori} as a pseudomonoid in $\cat{DblCat}$, and
define:
\begin{Defn}
A \defn{monoidal double morphism} between monoidal double
categories $\db K$ and $\db L$ is a (lax) map of pseudomonoids
$\db K \to \db L$ in $\cat{DblCat}$.
\end{Defn}
\noindent Again, the following is entirely routine:
\begin{Prop}\label{dblmorphalt}
Giving a monoidal double morphism $F \colon \db K \to \db L$ is
equivalent to giving a double morphism $F \colon \db K \to \db L$
such that
\begin{itemize*}
\item $F_0$ and $F_1$ are lax monoidal functors;
\item The equalities $sF_1 = F_0s$ and $tF_1 = F_0t$ hold as equalities of lax monoidal functors;
\item The natural transformations \[\m \colon F_1(\thg) \otimes F_1(?) \to F_1(\thg \,\otimes\, ?) \quad \text{and} \quad \e \colon \I_{F_0(\thg)} \to F_1(\I_{(\thg)})\]
are lax monoidal natural transformations (where we observe that
all the functors in question are indeed lax monoidal functors; for
instance, $F_1(\thg) \otimes F_1(?)$ is the composite
\[K_1 \mathbin{\vphantom{\times}_s\mathord{\times}_t} K_1 \xrightarrow{F_1 \mathbin{\vphantom{\times}_s\mathord{\times}_t} F_1}
L_1 \mathbin{\vphantom{\times}_s\mathord{\times}_t} L_1
\xrightarrow{\otimes} L_1
\]
which is the composite of a lax monoidal and a strong monoidal
functor as required).
\end{itemize*}
\end{Prop}

\noindent We can now define notions of
\emph{monoidal double homomorphism}, \emph{opmonoidal double
  morphism}, \emph{opmonoidal double opmorphism}, and so on.

Let us also note the correct notion of vertical transformation between monoidal double morphisms:
\begin{Defn}
A \defn{monoidal vertical transformation} between monoidal double
morphisms $F$, $G \colon \db K \to \db L$ is a pseudomonoid
transformation $F \Rightarrow G$ in $\cat{DblCat}$.
\end{Defn}
\begin{Prop}
Giving a monoidal vertical transformation $\alpha \colon F \Rightarrow G$ is equivalent
to giving a vertical transformation $\alpha \colon F \Rightarrow G$ such that $\alpha_0$ and $\alpha_1$
are monoidal transformations.
\end{Prop}

\noindent Straightforwardly, monoidal double categories, monoidal
double morphisms and monoidal vertical transformations form a
2-category $\cat{MonDblCat}$, along with all the expected variants:
$\cat{MonDblCat}_\psi$, $\cat{OpMonDblCat}$, $\cat{OpMonDblCat}_o$,
and so on.

\subsection{The monoidal double category $[\db K, \db K]_\psi$}
Given a small category $\cat{C}$, the endofunctor category
$[\cat{C}, \cat{C}]$ acquires the structure of a monoidal category.
We shall see in this section that a similar result holds for pseudo
double categories, namely, that the endohom double category $[\db K,
\db K]_\psi$ is naturally a monoidal double category.

Just as with transformations between morphisms of bicategories, there are
two canonical choices for the composite of two horizontal transformations
\[\A \colon A_s \Tor A_t \colon \db K \to \db K \quad \text{and}
\quad \B \colon B_s \Tor B_t \colon \db K \to \db K\text,\] namely
\[A_t \B \otimes \A
B_s \quad \text{and} \quad \A B_s \otimes A_t \B.\] As with the
bicategorical case, it makes no material difference which we choose:
\begin{Prop}
There are canonical invertible special modifications
\[i_{\A, \B} \colon A_t \B \otimes \A
B_s \Rrightarrow \A B_s \otimes A_t\B\text,\] natural in $\A$ and
$\B$.
\end{Prop}
\begin{proof}
We take $i_{\A, \B}$ to have central natural transformation
$A_{B_c(\thg)}$; so the component of $i_{\A, \B}$ at $X$ is given by
\[A_{\B X} \colon A_t\B X \otimes \A B_sX \to \A B_sX \otimes A_t\B X\text.\]
Visibly this is compatible with source and target, whilst the other
modification axiom is a long diagram chase using the axioms for $\A$ and $\B$.
For the naturality of these maps in $\A$ and $\B$, suppose we are given
modifications $\b \alpha \colon \A \Rrightarrow \C$ and $\b \beta \colon \B
\Rrightarrow \D$. Then we require the following diagrams to commute for all $X
\in K_0$:
\[\cd[@C+2.5em@R+1em]{
 A_t\B X \otimes \A B_sX
  \ar[d]_{A_{\B X}}
  \ar[r]^-{A_t\b \beta_X \otimes \A(\beta_s)_X} &
 A_t\D X \otimes \A D_sX
  \ar[d]^{A_{\D X}}
  \ar[r]^-{(\alpha_t)_{\D X} \otimes \b \alpha_{D_sX}} &
 C_t\D X \otimes \C D_sX
  \ar[d]^{C_{\D X}} \\
 \A B_sX \otimes A_t\B X
  \ar[r]_-{\A(\beta_s)_X \otimes A_t\b \beta_X} &
 \A D_sX \otimes A_t\D X
  \ar[r]_-{\b \alpha_{D_sX} \otimes (\alpha_t)_{\D X}} &
 \C D_sX \otimes C_t\D X\text.
}
\]
But the left-hand square is a naturality square for $A_{(\thg)}$
whilst the right-hand square is one of the axioms for $\b \alpha$; and
hence we are done.
\end{proof}

\begin{Prop}
The double category $[\db K, \db K]_\psi$ is a monoidal double
category.
\end{Prop}
\begin{proof}\hfill
\begin{itemize}
\item \textbf{Monoidal structure on $[\db K, \db K]_{v \psi}$}: Observe that
  this is the hom-category $\cat{DblCat}_\psi(\db K, \db K)$ in the
  2-category $\cat{DblCat}_\psi$, and hence is equipped with a strict
  monoidal structure;
\item \textbf{Monoidal structure on $[\db K, \db K]_{h \psi}$}: We take for
  the tensor unit $\b\unit$, the object
\[\b\unit = \I_{\id} \colon \id \Tor \id\text.\]
The tensor product is given as follows:
\begin{itemize}
\item \textbf{On objects}: given $\A \colon A_s \Tor A_t$ and $\B \colon B_s \Tor B_t$,
we take
\[\A \bten \B = \A B_t \otimes A_s \B \colon A_sB_s \Tor A_tB_t\text.\]
Explicitly, this has components
\[
  (\A \ten \B)(X)  = \cdl{A_sB_sX \ar[r]|-{\object@{|}}^-{A_s\B X} & A_sB_tX \ar[r]|-{\object@{|}}^-{\A B_tX} & A_tB_tX\text.}
\]

\item \textbf{On maps}:
Given
$\b \alpha \colon \b A \Rrightarrow \b C$ and $\b \beta \colon \b B \Rrightarrow \b D$, we take
\[\b \alpha \bten \b \beta = \b \alpha \beta_t \otimes \alpha_s \b \beta \colon \A B_t \otimes A_s \B \Rrightarrow \C D_t \otimes C_s \D\text.\]
\end{itemize}
The functoriality of $\bten$ is immediate from the functoriality of $\otimes$
and of the whiskering operations. We must now exhibit the unitality and
associativity coherence constraints in $[\db K, \db K]_{v\psi}$. For unitality,
we have that $\b\unit \bten \A = \I_\id A_t \otimes \A$ and $\A \bten \b\unit =
\A \otimes A_s \I_{\id}$, and so we give $\rho_\A$ and $\lambda_\A$ by the
special invertible modifications
\begin{align*}
& \cd{
 \A
  \ar@3[r]^-{\l_\A} &
 \I_{A_t} \otimes \A
  \ar@3[r]^-{\id} &
 \I_{\id}A_t \otimes \A
}\\
\text{and } & \cd{
 \A
  \ar@3[r]^-{\r_\A} &
 \A \otimes \I_{A_s}
  \ar@3[r]^-{\A \otimes \e_{\A}} &
 \A \otimes A_s \I_{\id}}
\end{align*} respectively. The naturality of these in $\A$ follows from the
naturality of $\l$, $\r$ and $\e$. For the associativity
modifications, suppose we are given $\A \colon A_s \Tor A_t$, $\B
\colon B_s \Tor B_t$ and $\C \colon C_s \Tor C_t$. Now we have
\begin{align*}
\A \bten (\B \bten \C) & = \A (B_t C_t) \otimes A_s(\B C_t \otimes
B_s \C)\\
\text{and } (\A \bten \B) \bten \C & = (\A B_t \otimes A_s \B) C_t \otimes (A_s
B_s) \C\text.
\end{align*}
Hence we take $\alpha_{\A, \B, \C}$ to be the special modification
\[\cd{
 \A (B_t C_t) \otimes A_s(\B C_t \otimes B_s \C)
  \ar@3[d]^{\ \A (B_t C_t) \otimes \m^{-1}_{\B C_t, B_s \C}} \\
 \A (B_t C_t) \otimes \big(A_s(\B C_t) \otimes A_s(B_s\C)\big)
  \ar@3[d]^{\ \id} \\
 (\A B_t) C_t \otimes \big((A_s \B) C_t \otimes (A_s B_s)\C\big)
  \ar@3[d]^{\ \a_{(\A B_t) C_t, (A_s \B) C_t, (A_s B_s)\C}} \\
 \big((\A B_t)C_t \otimes (A_s \B) C_t\big) \otimes (A_s B_s) \C
  \ar@3[d]^{\ \id} \\
 (\A B_t \otimes A_s \B) C_t \otimes (A_s B_s) \C\text.}\]
The naturality of these components in $\A$, $\B$ and $\C$ follows from the
naturality of $\m$ and $\a$; and a routine diagram chase using the coherence
axioms for $\l$, $\r$, $\a$, $\m$ and $\e$ shows that $\alpha$, $\rho$ and
$\lambda$ satisfy the associativity pentagon and the unit triangles.

\item \textbf{$s$ and $t \colon [\db K, \db K]_{h \psi} \to [\db K, \db K]_{v \psi}$ are strict monoidal}: this is immediate from above.

\item \textbf{$\I \colon [\db K, \db K]_{v \psi} \to [\db K, \db K]_{h \psi}$ is strong monoidal}: We observe that $\I_\unit = \b\unit$,
so that $\I$ is strict monoidal with respect to the unit. For the binary tensor
$\bten$, we have $\I_F \bten \I_G = \I_FG \otimes F\I_G$, and so we take $u_{F,
G} \colon \I_{FG} \Rrightarrow \I_F \bten \I_G$ to be the special invertible
modification
\[
\cd{
 \I_{FG}
  \ar@3[r]^-{\ \l_{\I_{FG}}} &
 \I_{FG} \otimes \I_{FG}
  \ar@3[r]^-{\ \id \otimes \e_G} &
 \I_FG \otimes F\I_G\text.
}\] Again, naturality in $F$ and $G$ follows from naturality of $\e$, and it's
easy to check that the three diagrams making $\I$ strong monoidal commute.

\item \textbf{$\otimes \colon [\db K, \db K]_{h \psi} \mathbin{\vphantom{\times}_s\mathord{\times}_{t}} [\db K, \db K]_{h \psi} \to [\db K, \db K]_{h \psi}$ is strong monoidal}:
Since $\I_\unit = \b\unit$, we can take
\[k_{\b\unit} \colon \b\unit \otimes \b\unit \to \b\unit\]
to be the canonical map $\r^{-1}_{\I_{\b \unit}} = \l^{-1}_{\I_{\b
\unit}}$. Now, suppose we are given horizontal transformations
\begin{align*}
\A \colon A_1 & \Tor A_2\text, & \A' \colon A_2 & \Tor A_3\text,\\
\B \colon B_1 & \Tor B_2& \text{and } \B' \colon B_2 & \Tor B_3\text.
\end{align*}
Then
\[(\A' \bten \B') \otimes (\A \bten \B) =
 (\A' B_3 \otimes A_2 \B') \otimes (\A B_2 \otimes A_1 \B)\]
whilst
\[(\A' \otimes \A) \bten (\B' \otimes \B) =
 (\A' \otimes \A)B_3 \otimes A_1(\B' \otimes \B)\text.\]
Therefore we take for $k_{\A', \B', \A, \B} \colon (\A' \bten \B')
\otimes (\A \bten \B) \to (\A' \otimes \A) \bten (\B' \otimes \B)$
the special invertible modification
\[\cd{
(\A' B_3 \otimes A_2 \B') \otimes (\A B_2 \otimes A_1 \B)
 \ar@3[d]^{\a} \\
\big(\A' B_3 \otimes (A_2 \B' \otimes \A B_2)\big) \otimes A_1 \B
 \ar@3[d]^{(\A' B_3 \otimes i_{\A, \B'}) \otimes A_1 \B} \\
\big(\A' B_3 \otimes (\A B_3 \otimes A_1 \B')\big) \otimes A_1 \B
 \ar@3[d]^{\a}s \\
(\A' B_3 \otimes \A B_3) \otimes (A_1 \B' \otimes A_1 \B)
 \ar@3[d]^{\id \otimes \m_{\B', \B}} \\
(\A' \otimes \A)B_3 \otimes A_1(\B' \otimes \B)
 }\]
where the maps labelled $\a$ are appropriate composites of associativity maps.
The naturality of the displayed map in all variables follows from the
naturality of $\a$, $i$ and $\m$. It's now a diagram chase to check that the
required coherence laws hold to make $\otimes$ strong monoidal.

\item \textbf{The natural transformations $\a$, $\l$ and $\r$ are strong monoidal transformations}: This is another routine diagram
chase.\qedhere
\end{itemize}
\end{proof}

\subsection{Monoidal comma double categories}
We now wish to mimic the result that, in the theory of clubs, tells us that
$[\cat C, \cat C] / S$ acquires a natural structure of monoidal category. As
there, we consider the lax limit of an arrow $\elt \X \colon 1 \to \db K$, but
this time in the 2-category $\cat{MonDblCat}$. Such an arrow amounts to a
\emph{monoidal monad} in $\cat{\db K}$:
\begin{Defn}
A \defn{monoidal monad} in the monoidal double category $\db K$
consists of:
\begin{itemize*}
\item A monad $(\X \colon X \tor X, \m, \e)$ in $\db K$;
\item Maps
\[\b \mu \colon \X \bten \X \to \X\text, \quad \b \eta \colon \b \unit \to \X\]
\[\mu \colon X \ten X \to X\text, \quad \text{and} \quad \eta \colon \unit \to X\]
\end{itemize*}
such that:
\begin{itemize*}
\item $s(\b \mu) = t(\b \mu) = \mu$ and $s(\b \eta) = t(\b \eta) = \eta$;
\item $(\X, \b \mu, \b \eta)$ is a monoid in the monoidal category $K_1$;
\item $(X, \mu, \eta)$ is a monoid in the monoidal category $K_0$;
\item The following diagrams commute:
\[\cd[@C+3em]{
  (\X \bten \X) \otimes (\X \bten \X)
    \ar[r]^-{k_{\X, \X, \X, \X}}
    \ar[d]_{\b \mu \otimes \b \mu} &
  (\X \otimes \X) \bten (\X \otimes \X)
    \ar[d]^-{\m \bten \m} \\
  \X \otimes \X
    \ar[d]_-{\m} &
  \X \bten \X
    \ar[d]^-{\b \mu} \\
  \X
    \ar[r]_-{\id} &
  \X
}\]
\[\cd{
  {\I_{X \ten X}}
     \ar[d]_{\I_{\mu}}
     \ar[r]^-{u_{X, X}} &
  {\I_X \bten \I_X}
     \ar[d]^{\e \bten \e} \\
  {\I_X}
     \ar[d]_{\e} &
  {\X \ten \X}
     \ar[d]^{\b \mu} \\
  \X
     \ar[r]_-{\id} &
  \X
}\quad \cd{
  {\b\unit \otimes \b\unit}
    \ar[r]^{k_{\b\unit}}
    \ar[d]_{\b \eta \otimes \b \eta} &
  {\b\unit}
    \ar[dd]^{\b \eta} \\
  {\X \otimes \X}
    \ar[d]_{\m} \\
  \X
    \ar[r]_{\id} &
  \X}
\quad \cd{
  {\b\I_\unit}
    \ar[r]^{u_{\unit}}
    \ar[d]_{\I_\eta} &
  {\b\unit}
    \ar[dd]^{\b \eta} \\
  {\I_X}
    \ar[d]_{\e} \\
  \X
    \ar[r]_{\id} &
  \X\text.}
\]
\end{itemize*}
\end{Defn}
\begin{Prop}\label{mdlcommadbl}
Let $\db K$ be a monoidal double category, and let $(\b X, \m, \e, \b \mu, \b \eta)$ be a monoidal monad in $\db K$.
Then the slice double category $\db K / \X$ can be equipped with the structure of a monoidal double category in
such a way as to become the lax limit of the arrow $\elt{\X} \colon 1 \to \db K$ in $\cat{MonDblCat}$.
\end{Prop}
\begin{proof}
We see that $\X$ and $X$ are monoids in the respective monoidal categories $K_1$ and $K_0$,
and therefore $K_1 / \X$ and $K_0 / X$ become monoidal categories.
It is straightforward to check that $s$ and $t$ are strict
monoidal with respect to this structure; for example, given $(\b U, \b
f)$ and $(\b U', \b f')$ in $K_1 / \X$, we have
$(\b U, \b f) \bten (\b U', \b f')$ given by
\[\cdl[@C+1em]{
  \b U \bten \b U' \ar[r]^-{\b f \bten \b f'} &
  \X \bten \X \ar[r]^-{\b \mu} &
  \X\text,
}
\]
whose image under $s$ is the object
\[\cdl[@C+2em]{
  U_s \ten U'_s \ar[r]^-{f_s \ten f'_s} &
  X \ten X \ar[r]^-{\mu} &
  X
}
\]
which is $(U_s, f_s) \ten (U'_s, f'_s)$ as required. It
remains to specify the invertible transformations $k$ and $u$ and
the invertible maps $k_{\b \unit}$ and $u_{\unit}$; the latter lift
straightforwardly from $\db K$, and the former we give as
follows:
\begin{align*}
 k_{(\b U, \b f), (\b U', \b f'), (\b V, \b g), (\b V', \b
 g')} &= k_{\b U, \b U', \b V, \b V'}\text, \\
 \text{and }
 u_{(U, f), (V, g)} &= u_{U, V}\text.
\end{align*}
That the required triangles commute for these to be maps in $K_1 / \X$ follows from the coherence
diagrams for $\X$; their naturality follows from the
naturality of $k$ and $u$ for $\db K$; and finally the
coherence diagrams that they are required to satisfy follow using
the coherence diagrams for $\X$ and $\db K$.
\end{proof}

\noindent In order to use this result in our theory of double clubs, we shall need the following:

\begin{Prop}
  Let $(S, \mu, \eta)$ be a double monad on a double category $\db K$.
  Then the monad $S\I$ in the monoidal double category $[\db K, \db
  K]_\psi$ is canonically a monoidal monad.
\end{Prop}
\begin{proof}
$S$ is a monad in $\cat{DblCat}_\psi$, and thus a monoid in
$\cat{DblCat}_\psi(\db K, \db K) = [\db K, \db K]_{v \psi}$. We equip the
object $S\I \in [\db K, \db K]_{h \psi}$ with monoid structure as follows.
Recall that $S\I$ is in fact the monad $S\I_{\id_{\db K}}$; so we give the unit
$\b \eta \colon \I_{\id_\db K} \Rrightarrow S\I$ by the modification
\[
\cd[@C+1.7em]{
 \b I_{\id_{\db K}}
   \ar@3[r]^-{\eta \I_{\id_{\db K}}} &
 S\I_{\id_{\db K}}\text.
}
\]
For the multiplication, observe first that we have $S\I \bten S\I =
(S\I_{\id_{\db K}})S \otimes S(S\I_{\id_{\db K}}) = S\I_S \otimes
S(S\I_{\id_{\db K}})$. Therefore we take for $\b \mu \colon S\I
\bten S\I \Rrightarrow S\I$ the modification
\[
\cd[@C+1.2em]{
 S\I_S \otimes S(S\I_{\id_{\db K}})
   \ar@3[r]^-{\m_{\I_S, S\I_{\id_{\db K}}}} &
 S(\I_S \otimes S\I_{\id_{\db K}})
   \ar@3[r]^-{S\l^{-1}_{S\I_{\id_{\db K}}}} &
 SS\I_{\id_{\db K}}
   \ar@3[r]^-{\mu \I_{\id_{\db K}}} &
 S\I_{\id_{\db K}}\text.
}
\]
It's straightforward to check that this makes $S\I$ into a monoid in
$[\db K, \db K]_{h \psi}$. Further, $s$ and $t$ send it to the monoid
$S$ in $[\db K, \db K]_{v \psi}$ as required. Finally, the diagrams
expressing the compatibility of the monoid and monad structure on $S$
are easily verified.
\end{proof}
\noindent Assembling the previous two results, we have:
\begin{Prop}
Given a double monad $(S, \eta, \mu)$ on a double category $\db
K$, the slice double category $[\db K, \db K]_\psi / S\I$ has a
natural structure of monoidal double category.
\end{Prop}

%% file: chapter5.tex
\section{Double clubs II}

We now have enough pseudo double category theory under our belt
to define the notion of a double club. First a few preliminaries:

\begin{Defn}Let $\db K$ and $\db L$ be double categories.
\begin{itemize}
\item  We say that $\db K$ is
  a \defn{vertically full sub-double category} of $\db L$ if there is a strict
  homomorphism $F \colon \db K \to \db L$ such that $F_0$ and
  $F_1$ exhibit $K_0$ and $K_1$ as full subcategories of $L_0$ and
  $L_1$.

\item  If $\db K$ and $\db L$ are monoidal double categories, we
say that
  $\db K$ is a \defn{sub-monoidal double category} of $\db L$ if there is a
  strict monoidal strict homomorphism $F \colon \db K \to \db L$ exhibiting
  $K_0$ and $K_1$ as subcategories of $L_0$ and $L_1$.
\end{itemize}
\end{Defn}
\noindent In particular, if $\db K$ is a vertically full
sub-double category of a monoidal double category $\db L$, then
$\db K$ can be made into a sub-monoidal double category of $\db L$
if and only the object sets of $K_0$ and $K_1$ are closed under
the binary and nullary tensors on $L_0$ and $L_1$ respectively.

\begin{Defn}
Let $(S, \eta, \mu)$ be a double monad on a double category $\db
K$. We say that $S$ is a \defn{double club} if:
\begin{itemize*}
\item $S$ has property (hps);
\item $\coll(S)$ is a sub-monoidal double category of $[\db K, \db K]_\psi / S\I$.
\end{itemize*}
\end{Defn}
\noindent Note that this is simply the natural generalisation of
Definition \ref{clubdef}: the extra requirement that condition
(hps) be satisfied is necessary to ensure that $\coll(S)$ exists
in the first place; in the plain category case, the existence of
the `category of collections' is automatic.

The above definition of a double club, though compact, is not very
easy to work with: but as with plain clubs, there is a more hands-on
description which greatly simplifies the task of applying the theory.

We begin by observing that if $(S, \eta, \mu)$ is a double monad on $\db K$,
then $(S_0, \eta_0, \mu_0)$ is a monad on $K_0$ and $(S_1, \eta_1, \mu_1)$ a
monad on $K_1$. Therefore it makes sense to ask whether or not $S_0$ and $S_1$
are clubs in the sense of Section \ref{chap:cart1} on their respective
categories, and once we have asked this, we may naturally ask whether this is
sufficient to make $S$ into a double club. In fact, as long as $S$ has property
(hps), the answer is yes:

\begin{Prop}\label{altclubdesc}
If $(S, \eta, \mu)$ is a double monad on $\db K$ such that:
\begin{itemize*}
\item $S$ has property (hps);
\item $S_0$ and $S_1$ are clubs on the categories $K_0$ and $K_1$ respectively,
\end{itemize*}
then $S$ is a double club.
\end{Prop}
\begin{proof}
We must check that $\coll(S)$ is a sub-monoidal double category of
$[\db K, \db K]_\psi / S\I$. Since $\coll(S)$ is a vertically full
sub-double category of $[\db K, \db K]_\psi / S\I$, it suffices to
check that:
\begin{itemize*}
\item $\coll(S)_0$ is closed under the monoidal structure on $[\db K, \db K]_{v \psi} / S$;
\item $\coll(S)_1$ is closed under the monoidal structure on $[\db K, \db K]_{h \psi} / S\I$.
\end{itemize*}
We begin with $\coll(S)_0$. We have evident forgetful functors
\[
\pi_i \colon [\db K, \db K]_{v \psi} / S \to [K_i, K_i] / S_i
\quad \text{(for $i = 0, 1$)}
\]
which are strict monoidal. Since $S_0$ and $S_1$ are clubs,
$Coll(S_i)$ is closed under the monoidal structure on $[K_i, K_i]
/ S_i$. But an object $A$ of $[\db K, \db K]_{v \psi}$ lies in
$\coll(S)_0$ just when its projections $\pi_i(A)$ lie in
$Coll(S_i)$; and hence we see that $\coll(S)_0$ is closed under
the monoidal structure on $[\db K, \db K]_{v \psi}$ as required.

Moving on to $\coll(S)_1$, we first show that the unit object $\b \eta \colon
\b I_{\id_{\db K}} \Rrightarrow S\I$ of $[\db K, \db K]_{h \psi}$ lies in
$\coll(S)_1$. By Proposition \ref{clubalt} and the fact that $S_0$ and $S_1$
are clubs, we have that $\eta_0$ and $\eta_1$ are cartesian natural
transformations; hence $\eta \colon \id_{\db K} \To S$ is a cartesian vertical
transformation. It remains to show that the central natural transformation of
$\b \eta$ is cartesian, i.e., that diagrams of the following form are
pullbacks:
\[
\cd[@C+2em]{
 {\I_X}
  \ar[r]^-{\I_!}
  \ar[d]_{\eta_{\I_X}} &
 {\I_1}
  \ar[d]^{\eta_{\I_1}} \\
 {S\I_{X}}
  \ar[r]_-{S\I_{!}} &
 {S\I_{1}}\text,
}
\]
which is just the cartesianness of $\eta_0$. We now show that $\coll(S)_1$ is
closed under the binary tensor product on $[\db K, \db K]_{h \psi}$. So suppose
we are given objects $(\b A, \b \alpha)$ and $(\b B, \b \beta)$ of
$\coll(S)_1$; then their tensor product is given by
\[
\cd[@+1em]{
A_sB_s \ar@2[d]_{\alpha_s\beta_s} \ar@2[r]^{\b A \bten \b B}|-{\object@{|}} \dthreecell{dr}{\b \alpha \bten \b \beta}& A_tB_t \ar@2[d]^{\alpha_t\beta_t}  \\
{SS} \ar@{=>}[d]_{\mu} \ar@{=>}[r]^{S\I \bten S\I}|-{\object@{|}} \dthreecell{dr}{\b \mu}& {SS} \ar@{=>}[d]^{\mu}  \\
S \ar@{=>}[r]|-{\object@{|}}_{S\I} & S\text, }\] so it suffices to
show that $\b \alpha \bten \b \beta$ and $\b \mu$ are cartesian
modifications. We begin with $\b \alpha \bten \b \beta$; the
cartesianness of $\alpha_s \beta_s$ and $\alpha_t \beta_t$ follows
from the fact that $S_1$ and $S_0$ are clubs on $K_1$ and $K_0$, and
so it suffices to check that the central natural transformation of
$\b \alpha \bten \b \beta$ is cartesian. This central natural
transformation has components
\[\cdl[@C+3em]{
  \A B_tX \otimes A_s\B X \ar[r]^{\b \alpha_{B_tX} \otimes (\alpha_s)_{\B X}} &
  S\I_{B_tX} \otimes S\B X \ar[r]^{S\I_{(\beta_t)_X} \otimes S\b \beta_X} &
  S\I_{SX} \otimes SS\I_{X}}\text.\]
So, consider the following diagram:
\[
\cd{
 \A B_tX
  \ar[d]_{\b \alpha_{B_tX}}
  \ar[r]^{\A B_t!} &
 \A B_t1
  \ar[d]^{\b \alpha_{B_t1}} \\
 S\I_{B_tX}
  \ar@/_4em/[ddd]_{S\I_{(\beta_t)_X}}
  \ar[d]_{S\e_X}
  \ar[r]_{S\I_{B_t!}} &
 S\I_{B_t1}
  \ar[d]^{S\e_1}
  \ar@/^4em/[ddd]^{S\I_{(\beta_t)_1}} \\
 SB_t\I_X
  \ar[d]_{S(\beta_t)_{\I_X}}
  \ar[r]_{SB_t\I_!} &
 SB_t\I_{1}
  \ar[d]^{S(\beta_t)_{\I_1}} \\
 SS\I_X
  \ar[d]_{S\e^{-1}_X}
  \ar[r]_{SS\I_!} &
 SS\I_{1}
  \ar[d]^{S\e^{-1}_1} \\
 S\I_{SX}
  \ar[r]_{S\I_{S!}} &
 S\I_{S1}\text.
}\] The top square is a pullback by cartesianness of $\alpha$, the
second and fourth are pullbacks since their vertical sides are
isomorphisms, and the third square is a pullback by cartesianness
of $\beta_t$ and because $S_1$ preserves cartesian natural
transformations into $S_1$. Therefore the outside edge of this
diagram is a pullback. Similarly, considering the diagram
\[
\cd{
 A_s\B X
  \ar[d]_{(\alpha_s)_{\B X}}
  \ar[r]^{A_s\B!} &
 A_s\B 1
  \ar[d]^{(\alpha_s)_{\B 1}} \\
 S\B X
  \ar[d]_{S\b \beta_X}
  \ar[r]_{S\B !} &
 S\B 1
  \ar[d]^{S\b \beta_1} \\
 SS\I_X
  \ar[r]_{SS\I_!} &
 SS\I_{1}\text,
}\] the top square is a pullback by cartesianness of $\alpha_s$, whilst the
bottom square is a pullback by cartesianness of $\beta$ and the fact that $S_1$
preserves cartesian transformations into $S_1$. Thus, forming the tensor
product of these two diagrams and applying condition (hps1), we see therefore
that the naturality squares for $(\b \alpha \ten \b \beta)_c$ are pullbacks as required.

Finally, we check that $\b \mu$ is a cartesian modification. By
Proposition \ref{clubalt} and the fact that $S_0$ and $S_1$ are
clubs, we have that $\mu_0$ and $\mu_1$ are cartesian natural
transformations; hence $\mu \colon SS \To S$ is a cartesian
vertical transformation. So we need only check that the central
natural transformation of $\b \mu$ is cartesian, for which we must
check that the outer edge of the following diagram is a pullback:
\[\cd[@C+2em]{
S\I_{SX} \otimes SS\I_X
 \ar[r]^-{S\I_{S!} \otimes SS!}
 \ar[d]_{\m_{\I_{SX}, S\I_X}} &
S\I_{S1} \otimes SS\I_1
 \ar[d]^{\m_{\I_{S1}, S\I_1}} \\
S(\I_{SX} \otimes S\I_X)
 \ar[r]^-{S(\I_{S!} \otimes S!)}
 \ar[d]_{S\l^{-1}_{S\I_X}} &
S(\I_{S1} \otimes S\I_1)
 \ar[d]^{S\l^{-1}_{S\I_1}} \\
SS\I_X
 \ar[r]^-{SS!}
 \ar[d]_{\mu_{\I_X}} &
SS\I_1
 \ar[d]^{\mu_{\I_1}} \\
S\I_X
 \ar[r]_-{S!} &
S\I_1\text.}\] Now, the bottom square is a pullback by
cartesianness of $\mu$, whilst all other squares are pullbacks
since they have isomorphisms along their vertical edges; hence the
outer edge is a pullback as required.
\end{proof}

%% file: chapter6.tex
\section{The double club for symmetric strict monoidal categories}
In Example \ref{etamu}, we saw that the monad on $\cat{Cat}$ for
symmetric strict monoidal categories extends to a double monad $S$
on $\dCat$. In Example \ref{clubsclubs}, we saw that this monad on
$\cat{Cat}$ is in fact a \emph{club} on $\cat{Cat}$. What we are now
in a position to show is that the double monad $S$ on $\dCat$ is
likewise a double club on $\dCat$.

Using Proposition \ref{altclubdesc}, this task is reduced to the
following: firstly, checking that $S_0$ and $S_1$ are clubs on their
respective categories, and secondly, showing that $S$ has property
(hps).  We have already seen in Example \ref{clubsclubs} that $S_0$
is a club on $\dCat_0 = \cat{Cat}$, and the following is a
straightforward but tedious calculation from the definitions:
\begin{Prop}\label{Sclub2}
The monad $(S_1, \eta_1, \mu_1)$ is a club on $\dCat_1$.
\end{Prop}
\noindent Therefore it remains only to show that $S$ satisfies property (hps), for which we shall
use the following two propositions:
\begin{Prop}\label{prop1}
Suppose that
\[\cd{
D \ar[r]^{j} \ar[d]_{k} & C \ar[d]^{g} \\
B \ar[r]_{f} & A}
\]
is a pullback in $\dCat_0$; then so is
\[\cd{
{\I_D} \ar[r]^{\I_j} \ar[d]_{\I_k} & {\I_C} \ar[d]^{\I_g} \\
\I_B \ar[r]_{\I_f} & {\I_A}}
\]
in $\dCat_1$.
\end{Prop}
\begin{proof}
  Viewing $\dCat_1$ as $\Cat / \b 2$, we see that the functor $\I_{(\
    )} \colon \dCat_0 \to \dCat_1$ sends $D$ to $(D \times \b 2)
  \xrightarrow{\pi_2} \b 2$, and is thus right adjoint to the domain
  functor $\Cat / \b 2 \to \Cat$.  Thus $\I_{(\
    )}$ preserves small limits and so \emph{a fortiori} the result.
\end{proof}

\begin{Prop}\label{prop2}
Let $A$ be a small groupoidal category and suppose we are given
pullback diagrams
\[(23):= \cd{{\b D}_{23} \ar[r]^{\b j_{23}} \ar[d]_{\b k_{23}} & {\b
    C_{23}} \ar[d]^{\b g_{23}} \\
{\b B}_{23} \ar[r]_{\b f_{23}} & {\I_A}} \quad \text{and}\quad
(12):= \cd{{\b D}_{12} \ar[r]^{\b j_{12}} \ar[d]_{\b k_{12}}& {\b
C_{12}}
  \ar[d]^{\b g_{12}} \\
{\b B}_{12} \ar[r]_{\b f_{12}} & {\I_A}}
\]
in $\dCat_1$ with
\[
s(12) = \cd{D_1 \ar[r]^{j_1} \ar[d]_{k_1} &
 C_1 \ar[d]^{g_1} \\
B_1 \ar[r]_{f_1} & A\text,} \quad \phantom{\text{and}} \quad t
(12) = \cd{D_2 \ar[r]^{j_2} \ar[d]_{k_2} &
 C_2 \ar[d]^{g_2} \\
B_2 \ar[r]_{f_2} & A\text,}
\]
\[
s (23) = \cd{D_2 \ar[r]^{j_2} \ar[d]_{k_2} &
 C_2 \ar[d]^{g_2} \\
B_2 \ar[r]_{f_2} & A\text,} \quad \text{and} \quad t (23) =
\cd{D_3 \ar[r]^{j_3} \ar[d]_{k_3} &
 C_3 \ar[d]^{g_3} \\
B_3 \ar[r]_{f_3} & A\text.}
\]
Suppose further that the arrow $f_2 \colon B_2 \to A$ is a
fibration; then the diagram
\[(13):= \cd{{\b D}_{23} \otimes {\b D}_{12} \ar[r]^{\b j_{23}
    \otimes \b j_{12}}
  \ar[d]_{\b k_{23} \otimes \b k_{12}}&
  {\b C_{23}} \otimes {\b C}_{12} \ar[d]^{\b g_{23} \otimes \b g_{12}} \\
{\b B}_{23} \otimes {\b B}_{12} \ar[r]_{\b f_{23} \otimes \b
f_{12}} & {\I_A} \otimes {\I_A}}
\]
is also a pullback.
\end{Prop}
\begin{proof}
First some notation; we shall use $b_i$, $c_i$ and $d_i$ to denote
typical elements of $B_i$, $C_i$ and $D_i$ (for $i = 1, \dots,
3$), and similarly use $a_i$ to denote elements of $A$,
 with the convention that
\[k_i(d_i) = b_i\text, \quad j_i(d_i) = c_i\text, \quad \text{and} \quad f_i(b_i) = a_i = g_i(c_i)\text.\]
So now, let $\b E = (E_1, E_2, E)$ be the pullback
\[\cd[@+1em]{{\b E} \ar[r]^{\b j'}
  \ar[d]_{\b k'}&
  {\b C_{23}} \otimes {\b C}_{12} \ar[d]^{\b g_{23} \otimes \b g_{12}} \\
{\b B}_{23} \otimes {\b B}_{12} \ar[r]_{\b f_{23} \otimes \b
f_{12}} & {\I_A} \otimes {\I_A}\text.}
\]
The universal property of pullback  induces a canonical arrow \[\b
u = (u_1, u_2, u) \colon {\b D}_{23} \otimes {\b D}_{12} \to \b
E\] in $\dCat_1$. It suffices to show that this map is an
isomorphism. Observe first that $s(13) = s(12)$ and $t(13) =
t(23)$, and thus that these projections are pullback diagrams in
$\cat{Cat}$. Thus we may take it that $E_1 = D_1$ and $E_2 = D_3$,
and that $u_1 = \id_{D_1}$ and $u_2 = \id_{D_3}$. Thus we need
only concern ourselves with the 2-cell $u$; we shall exhibit an
inverse $v$ for this 2-cell. First, let us describe explicitly
what $u$ does. A typical element of $\b D_{23} \otimes \b
D_{12}(d_3; d_1)$ looks like
\[\big((\alpha, \gamma)\otimes(\beta, \delta)\big) = \big((b_3, c_3) \xtor{(\alpha, \gamma)} (b_2, c_2)\big) \otimes
\big((b_2, c_2) \xtor{(\beta, \delta)} (b_1, c_1)\big)\] where
$\alpha \colon b_3 \tor b_2$, $\beta \colon b_2 \tor b_1$, $\gamma
\colon c_3 \tor c_2$, and $\delta \colon c_2 \tor c_1$ satisfy
\[a_3 \xrightarrow{f_{23}(\alpha)} a_2 = a_3 \xrightarrow{g_{23}(\gamma)} a_2 \quad \text{and} \quad
 a_2 \xrightarrow{f_{12}(\beta)} a_1 = a_2 \xrightarrow{g_{23}(\delta)} a_1\text,\]
whilst a typical element $\big((\alpha \otimes \beta), (\gamma \otimes \delta)\big)$ of $\b E(d_3; d_1)$ looks like
\[\big(
(b_3 \xtor{\alpha} b) \otimes (b \xtor{\beta} b_1), (c_3
\xtor{\gamma} c) \otimes (c \xtor{\delta} c_1)\big)\] where
\[a_3 \xrightarrow{f_{23}(\alpha)} f_2(b)
\xrightarrow{f_{12}(\beta)} a_1 = a_3 \xrightarrow{g_{23}(\gamma)}
g_2(c) \xrightarrow{g_{12}(\delta)} a_1\] in $A$. Then the 2-cell
$u$ has components given by
\begin{align*}
u_{d_3, d_1} \colon \b D_{23} \otimes \b D_{12}(d_3; d_1) & \rightarrow \b E(d_3; d_1)\\
\big((\alpha, \gamma)\otimes(\beta, \delta)\big) & \mapsto (\alpha
\otimes \beta, \gamma \otimes \delta)\text.
\end{align*}
Now let us construct the promised inverse $v$ for this 2-cell.
Suppose we are given an element $(\alpha \otimes \beta, \gamma
\otimes \delta) \in \b E(d_3; d_1)$; we must send this to an
element of ${\b D_{23} \otimes \b D_{12}}(d_3; d_1)$. So consider
the map
\[\psi := f_2(b) \xrightarrow{f_{23}(\alpha)^{-1}} a_3
\xrightarrow{g_{23}(\gamma)} g_2(c)\] in $A$. The functor $f_2
\colon B_2 \to A$ is a fibration and $A$ is a groupoid; thus $f_2$
is also a \emph{co}fibration, and so we can lift the displayed map
to a cocartesian arrow $\hat \psi \colon b \to \psi^* b$ in $B_2$;
and since $\psi$ is invertible, so is $\hat \psi$. So now we set
$v\big((\alpha \otimes \beta, \gamma \otimes \delta)\big)$ to be
\[
\big((b_3, c_3) \xtor{(\hat \psi \alpha, \gamma)} (\psi^*b, c)\big)
\otimes \big((\psi^* b, c) \xtor{(\beta \hat \psi^{-1}, \delta)} (b_1,
c_1)\big)\text.\] For this to be well-defined we need to check firstly
that it does indeed map into $\b D_{23} \otimes \b D_{12}(d_3; d_1)$;
and secondly that it is independent of the choice of representative
for $(\alpha \otimes \beta, \gamma \otimes \delta)$, both of which are
fairly tedious calculations which we therefore omit.

We must also check that $v$ is indeed inverse to $u$. We have $u\big((\alpha,
\gamma) \otimes (\beta, \delta)\big) = (\alpha \otimes \beta,
\gamma \otimes \delta)$, and thus $v\big(u\big((\alpha,
\gamma) \otimes (\beta, \delta)\big)\big)$ is given by
\[\big((b_3, c_3)
\xtor{(\hat \psi \alpha, \gamma)} (\psi^*b_2, c_2)\big) \otimes
\big((\psi^* b_2, c_2) \xtor{(\beta \hat
  \psi^{-1}, \delta)} (b_1, c_1)\big)\text,\]
where $\psi := g_{23}(\gamma) \circ f_{23}(\alpha)^{-1}$. But by
definition of $\b D_{23} \otimes \b D_{12}$, we have
$f_{23}(\alpha) = g_{23}(\gamma) \colon a_3 \to a_2$, and thus
\[vu\big((\alpha,
\gamma) \otimes (\beta, \delta)\big) = \big((\alpha, \gamma)
\otimes (\beta, \delta)\big)\] as required. Conversely, given
$(\alpha \otimes \beta, \gamma \otimes \delta)$ in $\b E(d_3;
d_1)$, we have that
\begin{align*}
uv(\alpha \otimes \beta, \gamma \otimes \delta) &= \big((\hat
\psi^{-1} \alpha) \otimes (\beta \hat \psi), \gamma \otimes
\delta\big)\\
&= \big(\alpha \otimes \beta, \gamma \otimes \delta\big)
\end{align*}
as required.
\end{proof}

\begin{Cor}\label{hps}
The homomorphism $S$ satisfies property (hps).
\end{Cor}
\begin{proof}
Condition (hps2) follows trivially from Proposition \ref{prop1}.
For (hps1), suppose we are given horizontally composable pullbacks
\[
\cd{
{\b A} \ar[r]^{\b p_1} \ar[d]_{\b p_2} & {\b B} \ar[d]^{\b f} \\
{S\b C} \ar[r]_{S!} & {S \I_1} } \quad \text{and} \quad
\cd{
{\b A'} \ar[r]^{\b p'_1} \ar[d]_{\b p'_2} & {\b B'} \ar[d]^{\b f'} \\
{S\b C'} \ar[r]_{S!} & {S \I_1}\text, }
\]
in $\dCat_1$. Then consider the diagram
\[\cd{{\b A'} \otimes {\b A} \ar[r]^{\b p'_1 \otimes \b p_1}
  \ar[d]_{\b p'_2 \otimes \b p_2}&
  {\b B'} \otimes {\b B} \ar[d]^{\b f' \otimes \b f} \\
{S \b C'} \otimes {S \b C} \ar[r]_{S ! \otimes S !} & {S\I_1}
\otimes {S\I_1} }
\]
We observe that $S1$ is a groupoid in $\cat{Cat}$, and that the arrow
$S! \colon SC_t \to S1$ in $\cat{Cat}$ is a fibration. We have an
isomorphism $S \I_1 \cong \I_{S1}$, and so can replace the
bottom-right vertex with $\I_{S1} \otimes \I_{S1}$; we now apply
Proposition \ref{prop2} to see that this square a pullback as
required.
\end{proof}

\begin{Cor}
The double monad $(S, \eta, \mu)$ is a double club on $\dCat$.
\end{Cor}
\begin{proof}
By Proposition \ref{hps}, $S$ has property (hps); and by
Proposition \ref{Sclub2}, $S_0$ and $S_1$ are
clubs on their respective categories. Therefore, by Proposition
\ref{altclubdesc}, $(S, \eta, \mu)$ is a double club on $\dCat$.
\end{proof}

In Example \ref{clubsclubs}, we also considered the clubs on
$\cat{Cat}$ for non-symmetric monoidal categories, and for
categories with finite products; in Example \ref{etamu}, we remarked
that they extended to double monads on $\dCat$. We leave it as an
exercise to the reader to show that these double monads are in fact
double clubs.


%% file: appendix.tex
\section*{Appendix A: Double equivalences}
We aim in this section to give an elementary characterisation of
equivalences in $\cat{DblCat}_\psi$. In fact, for very little extra
effort, we can garner significant extra generality by giving a
characterisation of adjunctions in $\cat{DblCat}$. A well-known
result in the theory of monoidal categories \cite{kelly:adj} says
that to give an adjunction in $\cat{MonCat}$, the 2-category of
monoidal categories, lax monoidal functors and monoidal
transformations, is to give an adjunction between the underlying
ordinary categories in $\cat{Cat}$ for which the left adjoint is
strong monoidal.

We shall produce a direct generalisation of this to pseudo double
categories, for which we need an analogue of `underlying ordinary
category'; more precisely, we need an appropriate analogue of the
2-category $\cat{Cat}$:

\begin{Defn}
We write $\cat{DblGph}$ for the 2-category
  $[\smash{\xymatrix@1{\bullet \ar@<2pt>[r] \ar@<-2pt>[r] & \bullet}}, \cat{Cat}]$.
\end{Defn}
There is an evident 2-functor $U \colon \cat{DblCat} \rightarrow
\cat{DblGph}$ which forgets horizontal structure, and so we may
speak of the `underlying double graph' of a double category.
\begin{Prop}
To give an adjunction $F \dashv G \colon \db L \to \db K$ in
$\cat{DblCat}$ is equivalent to giving an adjunction $F \dashv G
\colon U\db L \to U \db K$ in $\cat{DblGph}$ together with the
structure of a double homomorphism on $F$.
\end{Prop}
\noindent Let us spell out explicitly what the right hand side of
the above amounts to:
\begin{itemize*}
\item A double homomorphism $F \colon \db K \to \db L$;
\item A map of double graphs $G \colon \db L \to \db K$;
\item Adjunctions $F_0 \dashv G_0$ and $F_1 \dashv G_1$ with
  unit and counit $(\eta_0, \epsilon_0)$ and $(\eta_1, \epsilon_1)$ respectively,
\end{itemize*}
such that $s \epsilon_1 = \epsilon_0 s$, $t \epsilon_1 = \epsilon_0 t$, $s
\eta_1 = \eta_0 s$ and $t \eta_1 = \eta_0 t$.
\begin{proof}
  On an abstract level, this proof runs as follows:
  the 2-functor $U \colon \cat{DblCat} \to \cat{DblGph}$ has a
  left 2-adjoint $F$, which gives the `free double category' on a
  given double graph. Now, the 2-category of strict algebras and strict
  algebra maps for the induced monad $UF$ on $\cat{DblGph}$ is precisely
  the 2-category of \emph{strict} double
  categories, whilst the 2-category of \emph{pseudo}-algebras and
  lax algebra maps is \emph{almost} the 2-category
  $\cat{DblCat}$; more precisely, it is the 2-category of
  `unbiased' (in the sense of \cite{leinster:higheroperads}) pseudo double categories, which come equipped with $n$-ary
  horizontal composition functors for all $n$. As in
  the bicategorical case, it is not too hard to show that this
  notion is essentially equivalent to the `biased' notion of
  pseudo double category that we have adopted.

  Now, the 2-category $\cat{DblGph}$ is complete and cocomplete as
  a 2-category, and hence by Section 6.4 of \cite{bkp:monads}, there is a
  2-monad $T'$ on $\cat{DblGph}$ whose \emph{strict} algebras are
  precisely the \emph{pseudo} algebras for the composite monad
  $T = UF$. Thus, we have a 2-monad $T'$ on $\cat{DblGph}$ whose
  category of strict algebras and lax algebra maps can be
  identified with $\cat{DblCat}$.

  But now we are in a position to apply Kelly's `doctrinal adjunction'; by Theorem 1.5 of \cite{kelly:adj},
  to give an adjunction in $\cat{DblCat}$ is precisely to give an adjunction between the underlying
  objects of $\cat{DblGph}$ for which the left adjoint is a pseudo map of
  $T'$-algebras; and to give such a map is essentially the same
  thing as giving a homomorphism of pseudo double categories.

  Now, there are many details missing from the above, and rather than attempt to fill them in,
  it will be easier to give a direct proof following \cite{kelly:adj}.
  So, suppose first we are given an adjunction $UF \dashv UG$ in $\cat{DblGph}$ for which the left adjoint
  is a double homomorphism; then it suffices to
  equip $G$ with comparison transformations $\m$ and $\e$,
  and to show that $\eta = (\eta_0, \eta_1)$ and $\epsilon =
  (\epsilon_0, \epsilon_1)$ become vertical transformations with respect to
  this data.
  So, suppose that $F$ has comparison transformations
\[
\cd[@C+2em]{
  K_1 \cotimes K_1
   \dtwocell{dr}{\m}
   \ar[d]_{\otimes}
   \ar[r]^{F_1 \cotimes F_1} &
  L_1 \cotimes L_1
   \ar[d]^{\otimes} \\
  K_1
   \ar[r]_{F_1}&
  L_1
} \quad \text{and} \quad \cd[@C+2em]{
  K_0
   \dtwocell{dr}{\e}
   \ar[d]_{\I}
   \ar[r]^{F_0} &
  L_0
   \ar[d]^{\I} \\
  K_1
   \ar[r]_{F_1}&
  L_1
}
\]
Then we give the comparison transformations for $G$ as the mates
\[
\cd[@C+2em]{
  L_1 \cotimes L_1
   \dtwocell{dr}{\overline{\m^{-1}}}
   \ar[d]_{\otimes}
   \ar[r]^{G_1 \cotimes G_1} &
  K_1 \cotimes K_1
   \ar[d]^{\otimes} \\
  L_1
   \ar[r]_{G_1}&
  K_1
} \quad \text{and} \quad \cd[@C+2em]{
  L_0
   \dtwocell{dr}{\overline{\e^{-1}}}
   \ar[d]_{\I}
   \ar[r]^{G_0} &
  K_0
   \ar[d]^{\I} \\
  L_1
   \ar[r]_{G_1} &
  K_1
}
\]
of $\m^{-1}$ and $\e^{-1}$ under the adjunctions $F_0 \dashv G_0$,
$F_1 \dashv G_1$ and $F_1 \cotimes F_1 \dashv G_1 \cotimes G_1$.
Explicitly, the components of these transformations at $(\X, \Y)$
and $X$ respectively are given as follows:-
\[
\cd{
  G\X \otimes G\Y \ar[d]^{\eta_{G\X \otimes G\Y}} \\
  GF(G \X \otimes G\Y) \ar[d]^{G\m^{-1}_{G\X, G\Y}} \\
  G(FG \X \otimes FG \Y) \ar[d]^{G(\epsilon_\X \otimes \epsilon_\Y)} \\
  G(\X \otimes \Y)
} \quad \text{and} \quad \cd{
  \I_{GX} \ar[d]^{\eta_{\I_{GX}}} \\
  GF(\I_{GX}) \ar[d]^{G\e^{-1}_{X}} \\
  GFG\I_X \ar[d]^{G\epsilon_{\I_X}}\\
  G\I_X\text.
}
\]
That this data is coherent follows automatically from the coherence axioms for
$F$ and the functoriality of mates, and it's now a straightforward exercise in
the calculus of mates, following \cite{kelly:adj}, to show that $\eta =
(\eta_0, \eta_1)$ and $\epsilon = (\epsilon_0, \epsilon_1)$ become vertical
transformations with respect to this data. Thus we have an adjunction in
$\cat{DblCat}$ as required.

Conversely, any adjunction $(F, G, \eta, \epsilon)$ in
$\cat{DblCat}$ gives rise to the data specified above; we need
only check that $F$ is a homomorphism, i.e., that its special
comparison maps are invertible. Suppose that the comparison maps
for $G$ are $\m'$ and $\e'$; then it's easy to check that their
mates $\overline{\m'}$ and $\overline{\e'}$ furnish us with
inverses for $\m'$ and $\e'$ (explicitly, these inverses are given
by
\[
\cd{
  F(\X \otimes \Y) \ar[d]^{F(\eta_\X \otimes \eta_\Y)} \\
  F(GF\X \otimes GF\Y) \ar[d]^{F\m'_{F\X, F\Y}} \\
  FG(F\X \otimes F\Y) \ar[d]^{\epsilon_{F\X \otimes F\Y}} \\
  F\X \otimes F\Y
} \quad \text{and} \quad \cd{
  F\I_X \ar[d]^{F\I_{\eta_X}} \\
  F\I_{GFX} \ar[d]^{F\e'_X} \\
  FG\I_{FX} \ar[d]^{\epsilon_{\I_{FX}}}\\
  \I_{FX}\text.\ \text)
}
\]
The only thing remaining to check is that these two processes are
mutually inverse. Suppose we are given an adjunction $(F, G, \eta,
\epsilon)$ in $\cat{DblCat}$; then we must show that we can
reconstruct this adjunction from the underlying adjunction in
$\cat{DblGph}$ together with the data for $F$.

This amounts to checking that the special comparison maps we
produce for $G$ are the ones we started with; but this is
immediate, since we take them to be $\overline{\m^{-1}}$ and
$\overline{\e^{-1}}$, which are $\overline{\overline{\m'}} = \m'$
and $\overline{\overline{\e'}} = \e'$ as required.
\end{proof}

\begin{Cor}\label{equivv}
Suppose we are given double categories $\db K$ and $\db L$, and:
\begin{itemize*}
\item A double homomorphism $F \colon \db K \to \db L$;
\item A map of double graphs $G \colon \db L \to \db K$
\end{itemize*}
together with natural isomorphisms $\eta_i \colon \id_{K_i} \cong G_iF_i$ and
$\epsilon_i \colon F_iG_i \cong \id_{K_i}$ ($i = 0, 1$), such that such that $s
\epsilon_1 = \epsilon_0 s$, $t \epsilon_1 = \epsilon_0 t$, $s \eta_1 = \eta_0
s$ and $t \eta_1 = \eta_0 t$. Then $\db K$ and $\db L$ are equivalent in
$\cat{DblCat}_\psi$.
\end{Cor}
\begin{proof}
To give this data is to give an equivalence in $\cat{DblGph}$, so
by replacing $\epsilon_1$ and $\epsilon_0$, we can make this into
an \emph{adjoint} equivalence in $\cat{DblGph}$. Now, applying the
previous result, we get an (adjoint) equivalence in
$\cat{DblCat}$; but now we note that the comparison special maps
for $G$ will be invertible, since they are constructed from a
composite of invertible maps, and hence that our equivalence is an
equivalence in $\cat{DblCat}_\psi$ as well.
\end{proof}

\section*{Appendix B: Whiskering and double clubs}

We have defined the concept of double club in terms of closure under
the structure of monoidal double category. However, we may also ask
about closure under the `whiskering' operations of Section 2.
\emph{Prima facie}, this may appear to be a strictly stronger
requirement, but in fact it follows from the definition of double club
given above.

We begin with a preliminary general result on endohom double
categories.  We saw how to construct the monoidal structure on $[\db
K, \db K]_\psi$ using the whiskering operations $G(\thg)$ and
$(\thg)G$. We can also to a certain extent go in the other direction,
and derive something like the whiskering homomorphisms from the
monoidal structure on $[\db K, \db K]_\psi$.
Indeed, given a homomorphism $G \colon \db K \to \db K$, we obtain
homomorphisms
\begin{align*}
(\thg) \ten \I_G \colon [\db K, \db K]_\psi &\cong [\db K, \db K]_\psi \times 1 \xrightarrow{\id \times \elt{\I_G}}
 [\db K, \db K]_\psi \times [\db K, \db K]_\psi \xrightarrow{\ten}
 [\db K, \db K]_\psi\\
 \I_G \ten (\thg) \colon [\db K, \db K]_\psi &\cong 1 \times [\db K, \db K]_\psi \xrightarrow{\elt{\I_G} \times \id }
 [\db K, \db K]_\psi \times [\db K, \db K]_\psi \xrightarrow{\ten}
 [\db K, \db K]_\psi\text.\end{align*} And these homomorphisms approximate
 the operation of whiskering by $G$ in the following sense:
 \begin{Prop}\label{whiskerapprox}
 There are canonical invertible vertical transformations
 \[l_G \colon G(\thg) \Rightarrow \I_G \ten (\thg) \quad \text{and}
 \quad r_G \colon (\thg)G \Rightarrow (\thg) \ten \I_G\] which are
 natural in $G$.
 \end{Prop}
 \begin{proof}
 We have $\big(G(\thg)\big)_0 = \big(\I_G \ten (\thg)\big)_0$ and
 $\big((\thg)G)_0 = \big((\thg) \ten \I_G\big)_0$, so we can take
 $(l_G)_0$ and $(r_G)_0$ to be identity natural transformations.
 For $(l_G)_1$ and $(r_G)_1$, observe that we have
 \begin{align*}
 \big(\I_G \ten (\thg)\big)_1 = \I_G(\thg)_t \otimes G(\thg) &=
 \I_{G(\thg)_t} \otimes G(\thg)\\\text{and }
 \big((\thg) \ten \I_G \big)_1
 & = (\thg)G \otimes (\thg)_s\I_G\text.\end{align*} Therefore we take $(l_G)_1$
 to be the natural transformation
 \[(l_G)_1 = \cdl[@+0.5em]{G(\thg) \ar@2[r]^-{\l_{G(\thg)}} & \I_{G(\thg)_t} \otimes G(\thg)}\]
 and $(r_G)_1$ to be the natural transformation
 \[(r_G)_1 = \cdl[@+0.5em]{(\thg)G \ar@2[r]^-{\r_{(\thg)G}}&  (\thg)G \otimes \I_{(\thg)_sG} \ar@2[r]^-{\id \otimes \e_G} & (\thg)G \otimes (\thg)_s\I_{G}}\text.\]
It's now routine diagram chasing to check
that $l$ and $r$ satisfy all the required axioms for a vertical transformation, and that they are natural in $G$ as required.
 \end{proof}

\begin{Prop}\label{whiskerlift}
Let $S$ be a double club, and let $(A, \alpha)$ be an object of
$\coll(S)$. Then the whiskering homomorphisms
\[(\thg)A \colon [\db K, \db K]_\psi \to [\db K, \db K]_\psi \quad \text{and} \quad A(\thg) \colon [\db K, \db K]_\psi \to [\db K, \db K]_\psi\]
lift to homomorphisms
\[(\thg)(A, \alpha) \colon \coll(S) \to \coll(S) \quad \text{and} \quad (A, \alpha)(\thg) \colon \coll(S) \to \coll(S)\text.\]
\end{Prop}
\begin{proof}
We give the details for $(A, \alpha)(\thg)$, since $(\thg)(A,
\alpha)$ follows similarly. Following Proposition
\ref{whiskerapprox}, we have the homomorphism $\I_{(A, \alpha)} \ten
(\thg) \colon \coll(S) \to \coll(S)$; further we have the invertible
special vertical transformation
\[
 l_A \colon A(\thg) \Rightarrow \I_A \ten (\thg) \colon
 \db K \to \db K
\]
So we give $(A, \alpha)(\thg)$ as follows. Its component $\big((A,
\alpha)(\thg)\big)_0 \coll(S)_0 \to \coll(S)_0$ is simply
$\big(\I_{(A, \alpha)} \ten (\thg)\big)_0 = (A, \alpha) \ten
(\thg)$, whilst $\big((A, \alpha)(\thg)\big)_1 \colon \coll(S)_1 \to
\coll(S)_1$ is given as follows:
\begin{itemize}
\item \textbf{On objects}: given $(\B, \b \beta)$ in $\coll(S)_1$, we take
$(A, \alpha)(\B, \b \beta)$ to be the modification
\[
\cd[@+1em]{A\B \ar@3[r]^-{(l_A)_\B} & \I_A \bten \B
\ar@3[r]^-{\I_{\alpha} \bten \b \beta} &
\I_S \bten S\I \ar@3[r]^-{\e \bten S\I} &
S\I \bten S\I \ar@3[r]^-{\m} & S\I\text.}
\]
The first modification above is cartesian since it is invertible,
whilst the remaining composite is $\I_{(A, \alpha)} \bten (\b B, \b
\beta)$, and hence cartesian since $S$ is a double club; thus the
entire composite is cartesian as required.
\item \textbf{On maps}: given $\b \delta \colon (\B, \b \beta) \to (\C, \b
\gamma)$, we take $(A, \alpha)(\b \delta)$ to be given by
\[A \b \delta \colon (A, \alpha)(\B, \b \beta) \to (A, \alpha)(\C, \b \gamma)\text.\]
That this map is compatible with the projections down to $S\I$ is an
easy diagram chase.
\end{itemize}
It's immediate that these definitions are compatible with source and
target; it remains to give the comparison maps $\m$ and $\e$, for which we simply take
\begin{align*}\e_{(B, \beta)} = \e_B & \colon \I_{AB} \Rrightarrow A\I_{B}\\ \text{and }
\m_{(\B, \b \beta), (\B', \b \beta')} = \m_{\B, \B'} & \colon A\B
\otimes A\B' \to A(\B \otimes \B')\text.\end{align*} That these maps
are compatible with the projections down to $S\I$ is another
straightforward diagram chase, whilst the coherence axioms for $\m$
and $\e$ follows from those for $A(\thg)$ on $[\db K, \db K]_\psi$.
\end{proof}
\noindent For completeness, we also observe the following:
\begin{Prop}
Let $S$ be a double club, and let $\gamma \colon (A, \alpha) \to
(B, \beta)$ be a vertical arrow of $\coll(S)$. Then the whiskering
vertical transformations
\[(\thg)\gamma \colon (\thg)A \Rightarrow (\thg)B \quad \text{and} \quad \gamma(\thg) \colon A(\thg) \Rightarrow B(\thg)\]
lift to vertical transformations
\[(\thg)\gamma \colon (\thg)(A, \alpha) \Rightarrow (\thg)(B, \beta) \quad \text{and} \quad \gamma(\thg) \colon (A, \alpha)(\thg) \Rightarrow (B, \beta)(\thg)\text.\]
\end{Prop}
\noindent The proof is straightforward: one must simply show that
the components of $\gamma(\thg)$ and $(\thg)\gamma$ are compatible
with the projections down to $S\I$.
